\documentclass[a4paper,11pt]{article}

\usepackage[utf8]{inputenc}
\usepackage[T1]{fontenc}

\usepackage{lmodern}
\usepackage{a4wide}
\usepackage{amsmath,amssymb,amsfonts}
\usepackage{color}
\usepackage[ruled,vlined,lines numbered]{algorithm2e}
\usepackage{hyperref}
\usepackage{tikz}
\usetikzlibrary{decorations.pathreplacing}
\usetikzlibrary{matrix}
\usetikzlibrary{patterns}
\usepackage{graphicx}
\usepackage[geometry]{ifsym}
\usepackage[caption=false]{subfig}
\usepackage{multirow}

\DeclareMathOperator{\diag}{diag}

\title{Parallel Prony's method with multivariate matrix pencil approach and its numerical aspect}
\author{Nela Bosner}

\begin{document}
\maketitle

\begin{abstract}
Prony's method is a standard tool exploited for solving many imaging and data analysis problems that result in parameter identification in sparse exponential sums
$$f(k)=\sum_{j=1}^{T}c_{j}e^{-2\pi i\langle t_{j},k\rangle},\quad k\in \mathbb{Z}^{d},$$
where the parameters are pairwise different $\{ t_{j}\}_{j=1}^{M}\subset [0,1)^{d}$, and $\{ c_{j}\}_{j=1}^{M}\subset \mathbb{C}\setminus \{ 0\}$ are nonzero. The focus of our investigation is on a Prony's method variant based on a multivariate matrix pencil approach. The method constructs matrices $S_{1}$, \ldots , $S_{d}$ from the sampling values, and their simultaneous diagonalization yields the parameters $\{ t_{j}\}_{j=1}^{M}$. The parameters $\{ c_{j}\}_{j=1}^{M}$ are computed as the solution of an linear least squares problem, where the matrix of the problem is determined by $\{ t_{j}\}_{j=1}^{M}$. Since the method involves independent  generation and manipulation of certain number of matrices, there is intrinsic capacity for parallelization of the whole computation process on several levels. Hence, we propose parallel version of the Prony's method in order to increase its efficiency. The tasks concerning generation of matrices is divided among GPU's block of threads and CPU, where heavier load is put on the GPU. From the algorithmic point of view, the CPU is dedicated to the more complex tasks of computing SVD, eigendecomposition, and solution of the least squares problem, while the GPU is performing matrix--matrix multiplications and summations. With careful choice of algorithms solving the subtasks, the load between CPU and GPU is balanced. Besides the parallelization techniques, we are also concerned with some numerical issues, and we provide detailed numerical analysis of the method in case of noisy input data. Finally, we performed a set of numerical tests which confirm superior efficiency of the parallel algorithm and consistency of the forward error with the results of numerical analysis.
\end{abstract}

% The joint eigenbasis is obtained from the eigendecomposition of a single matrix that is random linear combination of  $S_{1}$, \ldots , $S_{d}$.
%On the first level, the tasks concerning generation of matrices is divided among GPU's block of threads and CPU, where heavier load is put on the GPU. On the second level, the individual threads are dealing with individual matrix elements. From the algorithmic point of view, the CPU is dedicated to the more complex tasks of computing SVD, eigendecomposition, and solution of the least squares problem, while the GPU is performing matrix--matrix multiplications and summations.

\section{Introduction}
Prony's method is a standard tool for parameter identification in sparse exponential sums
\begin{equation}\label{eq_exp_sum}
f(k)=\sum_{j=1}^{m}c_{j}e^{-2\pi \iota \langle t_{j},k\rangle},\quad k\in \mathbb{Z}^{d},
\end{equation}
for small $m$, occurring in many imaging and data analysis problems. The task is to determine pairwise different parameters $\{ t_{j}\}_{j=1}^{m}\subset [0,1)^{d}$, and nonzero coefficients $\{ c_{j}\}_{j=1}^{m}\subset \mathbb{C}\setminus \{ 0\}$ from relatively few sampling values. The focus of our investigation is on a Prony's method variant based on a multivariate matrix pencil approach presented in \cite{prony_matrix}. Our intention is to give a numerical aspect of this method, as well as to propose a more efficient way of its implementation.

%The method from \cite{prony_matrix} constructs matrices $S_{1}$, \ldots , $S_{d}$ from the sampling values, simultaneously diagonalizes them, and the parameters $\{ t_{j}\}_{j=1}^{m}$ are derived from their eigenvalues. The joint eigenbasis is obtained from the eigendecomposition of a single matrix that is random linear combination of $S_{1}$, \ldots , $S_{d}$. The coefficients $\{ c_{j}\}_{j=1}^{m}$ are computed as solution of an linear least squares problem, where $\{ t_{j}\}_{j=1}^{m}$ determine elements of the problem matrix.

Details of the algorithm are following.
For $n\in \mathbb{N}$, let $I_{n}=\{ 0,\ldots ,n\}^{d}$ with fixed ordering of the elements, let us build the $N\times N$ matrices
$$T=[f(k-h)]_{k,h \in I_{n}},\quad T_{\ell}=[f(k-h+e_{\ell})]_{k,h \in I_{n}},\ \ell =1,\ldots ,d$$
from sampling values, where $N=\# I_{n}=(n+1)^{d}$, and $e_{\ell}$ is the $\ell$-th unit vector. If $T$ has rank $m$, then one should compute reduced singular value decomposition (SVD)
\begin{equation}\label{eq_T_svd}
T=U\Sigma V^{*},
\end{equation}
where $\Sigma \in \mathbb{R}^{m\times m}$ is diagonal positive definite, $U,V\in \mathbb{C}^{N\times m}$ and $U^{*}U=V^{*}V=I\in \mathbb{C}^{m\times m}$. Let us define the set of $m\times m$ matrices
\begin{equation}\label{eq_generateSl}
S_{\ell}=U^{*}T_{\ell}V\Sigma^{-1},\quad \ell =1,\ldots, d.
\end{equation}
Theorem 2.1 from \cite{prony_matrix} states that if $n\ge \frac{K_{d}}{\min_{i\ne j}\| z_{i}-z_{j}\|}$, where $K_{d}$ is an absolute constant that only depends od $d$ and is specified in \cite{prony_multivar_gen}, then $T$ has rank equal to the number of terms $m$ in the exponential sum (\ref{eq_exp_sum}), and $S_{1}$,\ldots ,$S_{d}$  are simultaneously diagonalizable. Furthermore, let us define, as in \cite{prony_matrix}
$$z_{j}=e^{-2\pi \iota t_{j}}=[e^{-2\pi \iota t_{j}(1)}\ \ldots \ e^{-2\pi \iota t_{j}(d)}]^{T},\quad j=1,\ldots ,m,$$
where $t_{j}(\ell )$ denotes the $\ell$-th component of the vector $t_{j}$, and
$$z_{j}^{k}=e^{-2\pi \iota \langle t_{j},k\rangle},\quad j=1,\ldots ,m,\ k\in I_{n},$$
then any regular matrix $W$ that simultaneously diagonalizes $S_{1}$,\ldots ,$S_{d}$ yields a permutation $\tau$ on $\{ 1,\ldots ,m\}$ such that
\begin{equation}\label{eq_diagonalizeSl}
W^{-1}S_{\ell}W=\diag (\langle z_{\tau (1)},e_{\ell}\rangle ,\ldots ,\langle z_{\tau (m)},e_{\ell}\rangle ), \quad \ell =1,\ldots ,d.
\end{equation}
From the proof of this theorem it is also easy to see that for
$$A=[z_{j}^{k}]_{j=1,\ldots,m,\ k\in I_{n}}\in \mathbb{C}^{m\times N},\quad f=[f(k)]_{k\in I_{n}}\in \mathbb{C}^{N},\quad c=[c_{j}]_{j=1,\ldots ,m}\in \mathbb{C}^{m},$$
equation $f(k)=\sum_{j=1}^{m}c_{j}z_{j}^{k},$ $k\in I_{n}$, can be written in the matrix form as
\begin{equation}\label{eq_ls_probl}
f=A^{T}c.
\end{equation}
Theorem 2.1 in \cite{prony_matrix} represents a core of the Prony's algorithm. Simultaneous diagonalization of $S_{1}$, \ldots , $S_{d}$ is obtained by finding $W$ that diagonalizes a random linear combination
$$C_{\mu}=\sum_{\ell =1}^{d}\mu_{\ell}S_{\ell},\quad \mu \in \mathbb{C}^{d},$$
for random $\mu$ such that $\| \mu\|_{2}=1$. The system of equations (\ref{eq_ls_probl}) is solved as a linear least squares problem. Putting all this together results in Algorithm \ref{alg_prony}, as proposed in \cite{prony_matrix}.

\begin{algorithm}[ht]
    \KwIn{$f(k)$, $k\in I$}
    \KwOut{$t_{\tau (1)},\ldots ,t_{\tau (m)}$ and $c_{\tau (1)},\ldots ,c_{\tau (m)}$}

    \vspace*{0.3cm}
    Build the matrices $T$, $T_{\ell}$, $\ell =1,\ldots,d$, and the vector $f$\;
    Compute the reduced SVD of $T$\;
    Compute the matrices $S_{1},\ldots ,S_{d}$\;
    Choose random $\mu \in \mathbb{S}_{\mathbb{C}}^{d-1}$ and compute a matrix $W$ that diagonalizes $C_{\mu}$\;
    Use $W$ to simultaneously diagonalize $S_{1},\ldots ,S_{d}$ and reconstruct $z_{\tau (1)},\ldots ,z_{\tau (m)}$\;
    Compute $t_{\tau (j)}$ as the principal value of $\log (z_{\tau (j)})$, $j=1,\ldots ,m$\;
    Solve $\text{argmin}_{c}\| A^{T}c-f\|_{2}$ to recover $c_{\tau (1)},\ldots ,c_{\tau (m)}$\;
\caption{Prony's method with the multivariate matrix pencil approach}\label{alg_prony}
\end{algorithm}

From numerical point of view this algorithm involves several nontrivial issues, such as computing singular value decomposition and spectral decomposition (diagonalization), as well as establishing numerical rank of a matrix. The both decompositions are based on iterative methods and require respectable number of operations, which depends cubically on the matrix dimension. On the other hand, in the Prony's algorithm they are not in the same position, since dimensions of the matrices are different. Spectral decomposition is performed on a matrix with small dimension $m$, and does not pose a problem. The most time consuming task in the algorithm is full SVD of the matrix $T$, whose dimension is $N$. In case when parameters $t_{j}$ are close to each other, and particularly when $d>1$, $N$ can be quite large. Closely related with SVD is rank determination, which in the floating point arithmetics is not a straightforward task. Especially, when we introduce noise into the sampling values, we cannot expect that the matrix $T$ has rank $m$ any more, but there will be a significant drop in the singular values. It seams that in many applications of the Prony's method, full SVD is computed with all $N$ singular values, and then the singular values are examined. This is a wasteful procedure, since it spends a large number of operations for a small number of required singular values and vectors.

For computing reduced SVD more suitable are Lanczos bidiagonalization method \cite{lanczos_bidiag_golub_kahan}, \cite{lanczos_bidiag_paige}, and block power method for SVD \cite{svd_power_method}. The Lanczos method computes partial bidiagonalization of a matrix, and if it completes it produces left and right singular subspaces, and determines the rank. The power method computes left and right singular subspaces for fixed number of dominant singular values. In case when we know an overestimation of the rank, and if we combine the method with a rank revealing factorization on the projected matrix, we can obtain the reduced SVD.

On the other hand, as described earlier,  the Prony's method involves independent generation and manipulation of certain number of matrices, hence there is intrinsic capacity for parallelization of the whole computation process on several levels, using CPU and GPU.

The goal of this paper is to give a numerical aspect of the Prony's algorithm, it's numerical analysis in the case of the noisy input data, and to propose algorithmic enhancements as well as its parallel implementation. The parallel algorithm is going to exploit a hybrid CPU-GPU environment in order to produce the most efficient implementation.

%----> Ovdje staviti što će biti u kojem poglavlju.

\section{Choice of efficient SVD method}
Let us first make a short analysis of execution times for all steps in Algorithm \ref{alg_prony}. The following table displays times expressed  in seconds, for specific parameters $d$, $m$, $n$, and $N$.

{\scriptsize
\begin{center}
\begin{tabular}{|r|r|r|r|r|r|r|r|r|r|r|r|r|}
\hline
$d$ & $m$ & $n$ & $N$ & $t_{T,T_\ell}$ & $t_{SVD}$ & $t_{S_\ell}$ & $t_{C_\mu}$ & $t_{eig}$ & $t_{z,t}$ & $t_{A}$ & $t_{LS}$ & $t_{total}$\\
\hline
\hline
3 &  5 & 20 & 9261 & 48.12 & 166.55 & 0.1100 & 0.000040 & 0.000071 & 0.000038 & 0.01207 & 0.000968 & 214.80 \\
\hline
3 & 10 & 20 & 9261 & 47.30 & 232.13 & 0.1199 & 0.000046 & 0.000124 & 0.000089 & 0.02406 & 0.002209 & 279.57 \\
\hline
3 & 15 & 20 & 9261 & 47.34 & 232.52 & 0.1290 & 0.000046 & 0.000221 & 0.000064 & 0.03597 & 0.003891 & 280.04 \\
\hline
3 & 20 & 20 & 9261 & 47.93 & 211.00 & 0.1495 & 0.000046 & 0.000359 & 0.000135 & 0.04800 & 0.006585 & 259.13 \\
\hline
\end{tabular}
\end{center}
}
Notation of the time fractions is as follows:
\begin{itemize}
\item $t_{T,T_\ell}$ --- time spent on generating matrices $T$ and $T_{\ell}$, $\ell =1,\ldots , d$, and vector $f$
\item $t_{SVD}$ --- time spent on computing full SVD of $T$
\item $t_{S_\ell}$ --- time spent on generating matrices $S_{\ell}$
\item $t_{C_\mu}$ --- time spent on computing matrix $C_{\mu}$
\item $t_{eig}$ --- time spent on diagonalizing $C_{\mu}$
\item $t_{z,t}$ --- time spent on computing vectors $z_{j}$ and $t_{j}$
\item $t_{A}$ --- time spent on computing matrix $A$
\item $t_{LS}$ --- time spent on solving least squares problem $\text{argmin}_{c}\| A^{T}c-f\|_{2}$
\item $t_{total}$ --- total execution time
\end{itemize}
As we can see from the table, the dominant task is full SVD of the matrix $T$ which consumes about 80\% of the total execution time in the presented example for $d=3$ and $n=20$.

To avoid computation of full SVD of a large matrix, and to reduce time spent on finding reduced SVD the first step in producing an efficient Prony's algorithm is the right choice of the SVD algorithm. Instead of computing the full SVD, the algorithm should be able to determine rank of the matrix and to compute singular subspaces directly.

\subsection{Lanczos bidiagonalization method}
The Lanczos bidiagonalization method is suitable for rank determination in case when we have no information on rank. Nevertheless, this method have some numerical issues that we will refer later. As described in \cite[Subsection 9.3.3]{matrix_comp}, suppose $U^{*}AV=B$ represents the full bidiagonalization of $A\in \mathbb{C}^{m\times n}$ ($m\ge n$), with $U=[u_{1},\ldots ,u_{n}]\in \mathbb{C}^{m\times n}$ such that $U^{*}U=I_{n}$, $V=[v_{1},\ldots ,v_{n}]\in \mathbb{C}^{n\times n}$ such that $V^{*}V=I_{n}$, and
$$B=\left[ \begin{array}{ccccc}
\alpha_{1} & \beta_{2} & & \cdots & 0\\
0 & \alpha_{2} & \ddots & & \vdots \\
 & \ddots & \ddots & \ddots & \\
\vdots & & \ddots & \ddots & \beta_{n}\\
0 & \cdots & & 0 & \alpha_{n}
\end{array} \right] \in \mathbb{R}^{n\times n}.$$
By comparing columns in the equations $AV=UB$ and $A^{*}U=VB^{T}$ we obtain the method displayed in Algorithm \ref{alg_lanczos}.

\begin{algorithm}[htb]
  \caption{Lanczos bidiagonalization method}\label{alg_lanczos}
\KwIn{$A\in \mathbb{C}^{m\times n}$, $p_{1}\in \mathbb{C}^{n}$}
\KwOut{rank $r$, $U_{r}\in \mathbb{C}^{m\times r}$, $V_{r+1}\in \mathbb{C}^{n\times (r+1)}$, $B_{r,r+1} \in \mathbb{R}^{r\times (r+1)}$}

$p_{1}=$ given vector\;
$\beta_{1}=\| p_{1} \|_{2}$; $v_{1}=p_{1}/\beta_{1}$; $u_{0}=0$; $i=1$\;
\While{$\beta_{i} \ne 0$}
{
    $r_{i}=Av_{i}-\beta_{i}u_{i-1}$; $\alpha_{i}=\| r_{i}\|_{2}$\;
    \If{$\alpha_{i} = 0$}
    {
        exit the loop\;
    }
    $u_{i}=r_{i}/\alpha_{i}$; $i=i+1$\;
    $p_{i}=A^{*}u_{i-1}-\alpha_{i-1}v_{i-1}$; $\beta_{i}=\| p_{i}\|_{2}$\;
    \If{$\beta_{i} \ne 0$}
    {
        $v_{i}=p_{i}/\beta_{i}$\;
    }
}
$r=i-1$\;

\end{algorithm}

In the matrix form the algorithm can be rewritten as:
\begin{align}
AV_{i}&=U_{i}B_{i}\label{eq_lanc_rec1}\\
A^{*}U_{i}&=V_{i+1}B_{i,i+1}^{T}\label{eq_lanc_rec2}, &i=1,2,\ldots ,
\end{align}
where $B_{i}=B(1:i,1:i)$ is top $i\times i$ submatrix of $B$, $B_{i,i+1}=B(1:i,1:i+1)$ is top $i\times (i+1)$ submatrix of $B$, $U_{i}=[u_{1}\ \ldots \ u_{i}]$, and $V_{i}=[v_{1}\ \ldots \ v_{i}]$. If $\alpha_{i+1}=0$ then $\text{span}\{ Av_{1},\ldots ,Av_{i+1}\} \subset\{ u_{1},\ldots ,u_{i}\}$ which implies rank deficiency. In that case $B_{i+1}$ is singular with the last row equal to zero, hence (\ref{eq_lanc_rec1}) transforms into
$$AV_{i+1}=U_{i}B_{i,i+1},$$
which together with (\ref{eq_lanc_rec2}) shows that, like in \cite{lanczos_bidiag_paige}
$$A^{*}AV_{i+1}=V_{i+1}B_{i,i+1}^{T}B_{i,i+1},\qquad AA^{*}U_{i}=U_{i}B_{i,i+1}B_{i,i+1}^{T},$$
and $\sigma (B_{i,i+1}^{T}B_{i,i+1})\subset \sigma(A^{*}A)$ and $\sigma (B_{i,i+1}B_{i,i+1}^{T})\subset \sigma (AA^{*})$, where $\sigma (X)$ denotes specter of the matrix $X$. So, $\alpha_{i+1}=0$ implies that $i$ nontrivial singular values of $B_{i,i+1}$ are also singular values of $A$. Suppose that
$$B_{i,i+1}=\bar{U}_{B}[\bar{\Sigma}_{B}, 0] \bar{V}_{B}^{T},\quad \bar{U}_{B}\in \mathbb{R}^{i\times i},\ \bar{U}_{B}^{T}\bar{U}_{B}=I,\ \bar{V}_{B}\in \mathbb{R}^{(i+1)\times (i+1)},\ \bar{V}_{B}^{T}\bar{V}_{B}=I,$$
is SVD of $B_{i,i+1}$ where $\bar{\Sigma}_{B} \in \mathbb{R}^{i\times i}$ is diagonal, with positive diagonal elements. Then, diagonal elements of $\bar{\Sigma}_{B}$ are also singular values of $A$, columns of $U_{i}\bar{U}_{B}$ are left singular vectors of $A$, and columns of $V_{i+1}\bar{V}_{B}(:,1:i)$ are right singular vector of $A$. The vector $V_{i+1}\bar{V}_{B}(:,i+1)$ is in null space of $A$.

There is another stopping criteria for the Lanczos algorithm. When $\beta_{i+1}=0$, (\ref{eq_lanc_rec2}) reduces to
\begin{equation}\label{eq_beta0}
A^{*}U_{i}=V_{i}B_{i}^{T},
\end{equation}
which together with (\ref{eq_lanc_rec1}) implies that singular values of $B_{i}$ are also singular values of $A$. If
$$B_{i}=U_{B}\Sigma_{B} V_{B}^{T},\quad U_{B}\in \mathbb{R}^{i\times i},\ U_{B}^{T}U_{B}=I,\ V_{B}\in \mathbb{R}^{i\times i},\ V_{B}^{T}V_{B}=I,$$
is SVD of $B_{i}$ then columns of $U_{i}U_{B}$ are left singular vectors of $A$, and columns of $V_{i}V_{B}$ are right singular vector of $A$.

It is interesting to investigate when and how the bidiagonalization is going to stop. Similar to the conclusion of Paige in \cite{lanczos_bidiag_paige}, if $v_{1}\in \mathcal{R}(A^{*})$ then by line 8 in Algorithm \ref{alg_lanczos} it follows that $v_{i}\in \mathcal{R}(A^{*})$ for all $i$, meaning that all $v_{i}$ are orthogonal to $\mathcal{N}(A)$. If the method stopped with $\alpha_{i+1}=0$, then $AV_{i+1}\bar{V}_{B}(:,i+1)=0$ implying that there exist a linear combination of $v_{1}$,\ldots ,$v_{i+1}$ which is in $\mathcal{N}(A)$. This is contradiction, hence in this case the algorithm has to stop with $\beta_{i+1}=0$ for $i\le r$ where $r=\text{rank}(A)$. Early stopping with $i<r$ can only occur if $v_{1}$ is a linear combination of $i$ right singular vectors of $A$, belonging to nontrivial singular values. If $v_{1}\notin \mathcal{R}(A^{*})$ and the method stopped with $\beta_{i+1}=0$, then by (\ref{eq_beta0}) and nonsingularity of $B_{i}$ it follows that $V_{i}=A^{*}U_{i}B_{i}^{-T}$ and all $v_{j}\in \mathcal{R}(A^{*})$ which is contradiction. Hence in this case the algorithm has to stop with $\alpha_{i+1}=0$ for $i\le r$. Early stopping can occur if $v_{1}$ is a linear combination of $i<r$ right singular vectors of $A$ and a vector from the null space of $A$.

If we start the bidiagonalization algorithm with random vector $p_{1}$, there is a great chance that it would not belong to a space spanned by only a few right singular vectors and vectors from null space. But nevertheless, if early stopping occur there is a way to continue the process. First of all, we have to detect early stopping. In case of $\alpha_{i+1}=0$ we have to check that $\text{span}\{ u_{1},\ldots,u_{i}\}=\mathcal{R}(A)$, and respectively that it is orthogonal complement of $\mathcal{N}(A^{*})$. It suffices to take a random vector $y$, orthogonalize it against $u_{1},\ldots,u_{i}$ by $y=y-U_{i}U_{i}^{*}y$ and check whether $y\in \mathcal{N}(A^{*})$. If condition on $y$ is satisfied we can conclude that we are finished, otherwise we stopped too early and we can proceed the algorithm with $u_{i+1}=y/\| y\|_{2}$. In case of $\beta_{i+1}=0$ we have to check that $\text{span}\{ v_{1},\ldots,v_{i}\}=\mathcal{R}(A^{*})$, and respectively that it is orthogonal complement of $\mathcal{N}(A)$. Hence, we take a random vector $w$, orthogonalize it against $v_{1},\ldots,v_{i}$ by $w=w-V_{i}V_{i}^{*}w$ and check whether $w\in \mathcal{N}(A)$. If condition on $w$ is not satisfied we stopped too early and we can proceed the algorithm with $v_{i+1}=w/\| w\|_{2}$. This way we can be sure that the algorithm stops for $i=r$, and we can determine the rank of $A$ from dimensions of its bidiagonal factor.

When the Lanczos bidiagonalization method is implemented in the finite precision arithmetic on a computer, then there is several numerical issues that have to be taken under consideration. Let $u$ be the unit roundoff, and let $U_{i}$, $V_{i+1}$ and $B_{i,i+1}$ be matrices computed in the finite precision arithmetic by Algorithm \ref{alg_lanczos}, which stopped as described in previous paragraph. Orthogonality among the columns of $U_{i}$ and $V_{i+1}$ is gradually lost, and as a consequence a rank of $B_{i,i+1}$ can be larger of $r$, since set of its singular values will contain false multiple copies of singular values even if they are single, or ghost singular values that appear between singular values of $B_{i,i+1}$. So, implementation of the Algorithm \ref{alg_lanczos} should include full reorthogonalization, where each vector $u_{i}$ and $v_{i}$ should be explicitly orthogonalized against all previous vectors $u_{j}$ and $v_{j}$, respectively, for $j=1,\ldots ,i-1$. This adds much to the algorithms operation count, but in case when we do not expect the rank to be large it is not too expensive. In this case, by Theorem 5 in \cite{lanczos_bidiag_reorth}, we can expect that
$$\hat{U}_{i}A\hat{V}_{i+1}^{*}=B_{i,i+1}+E_{i,i+1},$$
where the columns of $\hat{U}_{i}$ form an orthonormal basis for $\text{span}\{ u_{1},\ldots,u_{i}\}$, the columns of $\hat{V}_{i+1}$ form an orthonormal basis for $\text{span}\{ v_{1},\ldots,v_{i+1}\}$, and the elements of $E_{i,i+1}$ are of order $\mathcal{O}(u\| A\|_{2})$. This means that stopping criteria has to be rephrased, instead for checking $\alpha_{i+1}=0$ and $\beta_{i+1}=0$ we will test for $\alpha_{i+1}\le \text{tol}$ and $\beta_{i+1}\le \text{tol}$. The tolerance is $\text{tol}=u\| A\|_{2}$, or, in case when we use the Lanczos method as a step of the Prony's method for noisy input data, we can increase the tolerance to the perturbation magnitude in order to catch the drop of singular values.

\subsection{Block power method for SVD}
In case when we can obtain a slight overestimation of the rank, we can use faster approach such as the block power method for SVD (see \cite{svd_power_method}). For a matrix $A\in \mathbb{C}^{m\times n}$ with rank $r$ the goal is to find orthonormal $\bar{U}\in \mathbb{C}^{m\times r}$ and $\bar{V}\in \mathbb{C}^{n\times r}$, and a nonsingular matrix $Q\in \mathbb{C}^{r\times r}$ such that
$$A=\bar{U}Q\bar{V}^{*}.$$
The reduced SVD of $A$ is then easily obtained from SVD of $Q=U_{Q}\Sigma V_{Q}^{*}$, so that columns of $U=\bar{U}U_{Q}$ represent left singular vectors, and columns of $V=\bar{V}V_{Q}$ represent right singular vectors of $A$. Our approach is to use power method with rank determination, displayed in Algorithm \ref{alg_power}.

\begin{algorithm}[htb]
\caption{Block power method for SVD}\label{alg_power}
    \KwIn{$A\in \mathbb{C}^{m\times n}$, overestimation of rank $r_{0}$, orthonormal $U_{0}\in \mathbb{C}^{m\times r_{0}}$, $V_{0}\in \mathbb{C}^{n\times r_{0}}$}
    \KwOut{rank $r$, $U\in \mathbb{C}^{m\times r}$, $V\in \mathbb{C}^{n\times r}$, $\Sigma \in \mathbb{R}^{r\times r}$}

    $Q_{0}=U_{0}^{*}AV_{0}$\;
    $R_{0}=AV_{0}-U_{0}Q_{0}$\;
    $k=0$\;
    \While{$\| R_{k}\|_{F}> \text{tol}\| A\|_{F}$}
    {
        $k=k+1$\;
        $\bar{U}_{k}=AV_{k-1}$\;
        Compute reduced QR factorization $\bar{U}_{k}=U_{k}R_{U_{k}}$\;
        $\bar{V}_{k}=A^{*}U_{k}$\;
        Compute reduced QR factorization with column pivoting $\bar{V}_{k}P_{V_{k}}=V_{k}R_{V_{k}}$\;
        Determine new rank $r_{k}$ from $R_{V_{k}}$\;
        $V_{k}=V_{k}(:,1:r_{k})$, $Q_{k}^{*}=(R_{V_{k}}P_{V_{k}}^{T})(1:r_{k},1:r_{k})$\;
        $R_{k}=AV_{k}-U_{k}Q_{k}$\;
    }
    Compute SVD $Q_{k}=U_{Q_{k}}\Sigma V_{Q_{k}}^{*}$\;
    $r=r_{k}$, $U=U_{k}U_{Q_{k}}$, $V=V_{k}V_{Q_{k}}$\;

\end{algorithm}

There is a slight difference between Algorithm \ref{alg_power} and the power method described in \cite{svd_power_method}. In lines 9 and 10 of Algorithm \ref{alg_power} there is QR factorization with column pivoting required for the rank determination instead of ordinary QR, and it has to be performed only once in the first step of the loop. Rank determination finds index $i$ such that $R_{V_{1}}(i:r_{0},i:r_{0})=0$, or in case of finite precision arithmetics $\|R_{V_{1}}(i:r_{0},i:r_{0})\|_{F}\le tol \| R_{V_{1}}\|_{F}$. Then we take $r_{1}=i-1$. Since we are searching for a sharp drop in singular values, meaning that $\sigma_{r_{1}}\gg \sigma_{r_{1}+1}$, we expect that the power method will converge only in a few steps.

Even if all columns of $V_{0}$ are linear combination of only $k<r$ right singular vectors and vectors from the null space of $A$, then the QR factorization in line 7 of Algorithm \ref{alg_power} will produce $U_{1}$ whose columns consist of $k$ left singular vectors and $r_{0}-k$ vectors which almost always have components in the rest $r-k$ left singular vectors. Then, we will detect that $r_{1}=r$ and the power method will produce expected result. In finite precision arithmetics there is no problem with orthogonality like in Lanczos method since QR factorization implemented with Householder reflectors is numerical stable, and we are going to use same tolerance $tol=u$, or higher in case of noisy data, for rank determination by QR factorization with pivoting.

Hence, the power method is a very convenient method for obtaining the reduced SVD in case when we know an overestimation of rank. Moreover, it employs mostly matrix-matrix multiplications which are the most efficient matrix operations, as well as QR factorizations which are implemented as efficient blocked algorithms, so we can expect that the SVD power method is the fastest SVD method especially for larger dimensions, and numerical tests are going to confirm that.

\section{Parallel algorithm}
The Prony's method involves independent generation and manipulation of certain number of matrices, so, as mentioned earlier, we can do that simultaneously on several levels, employing both CPU and GPU. On the first level of parallelization, the tasks concerning generation of matrices is divided among GPU's block of threads and CPU, where heavier load is put on the GPU. On the second level, the individual threads are dealing with individual matrix elements.

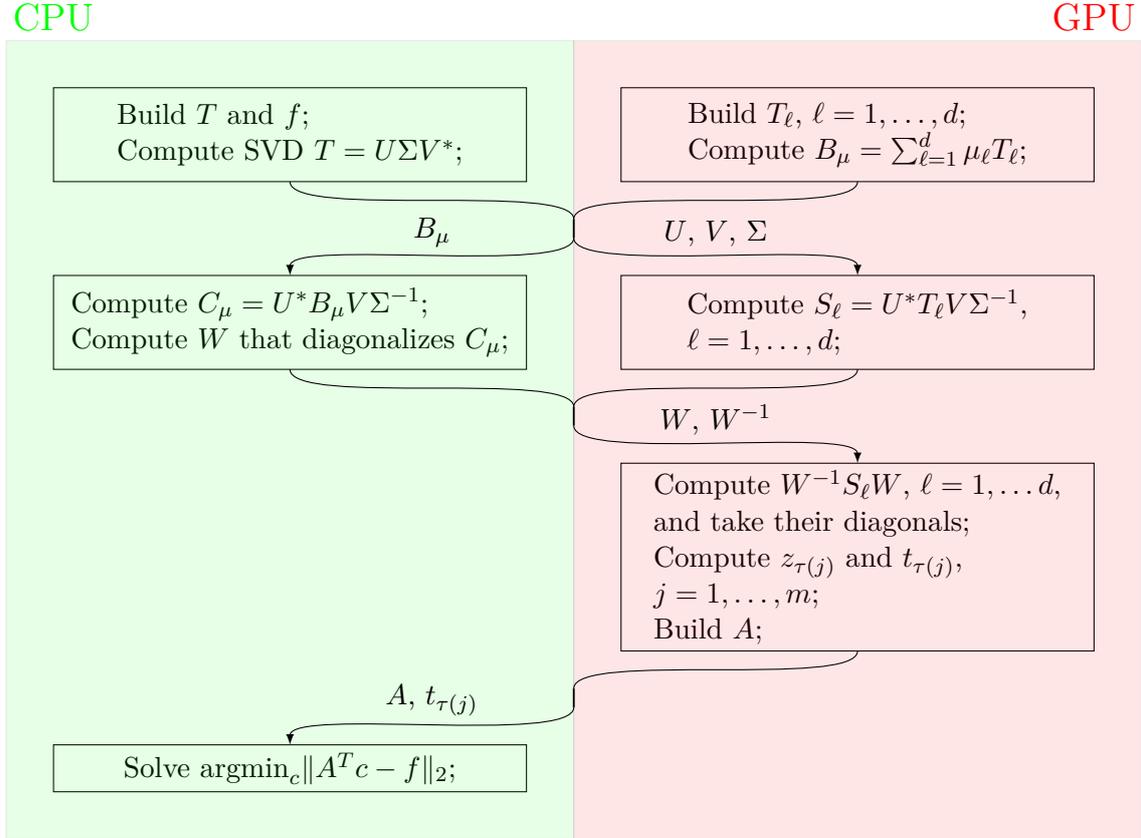
\begin{figure}[htb]
\begin{center}
\begin{tikzpicture}[x=\textwidth / 10,y=\textwidth / 10,xscale=0.8,yscale=0.8]
\draw (0,10) rectangle (5,9);
\draw (6,10) rectangle (11,9);
\draw (2.5,9.5) node[align=left] {Build $T$ and $f$;\\ Compute SVD $T=U\Sigma V^{*}$;};
\draw (8.5,9.5) node[align=left] {Build $T_{\ell}$, $\ell = 1,\ldots ,d$;\\ Compute $B_{\mu}=\sum_{\ell =1}^{d}\mu_{\ell}T_{\ell}$;};

\draw (0,8) rectangle (5,7);
\draw (6,8) rectangle (11,7);
\draw (2.5,7.5) node[align=left] {Compute $C_{\mu}=U^{*}B_{\mu}V\Sigma^{-1}$;\\ Compute $W$ that diagonalizes $C_{\mu}$;};
\draw (8.5,7.5) node[align=left] {Compute $S_{\ell}=U^{*}T_{\ell}V\Sigma^{-1}$,\\ $\ell =1,\ldots, d$;};

\draw(6,6) rectangle (11,4);
\draw (8.5,5) node[align=left] {Compute $W^{-1}S_{\ell}W$, ${\textstyle \ell =1,\ldots d}$,\\ and take their diagonals;\\ Compute $z_{\tau (j)}$ and $t_{\tau (j)}$,\\ $j=1,\ldots ,m$;\\ Build $A$;};

\draw (0,3) rectangle (5,2.5);
\draw (2.5,2.75) node[align=left] { Solve $\text{argmin}_{c}\| A^{T}c-f\|_{2}$;};

\filldraw[very nearly transparent, draw=black, fill=green] (-0.5,10.5) rectangle (5.5,2);
\filldraw[very nearly transparent, draw=black, fill=red] (5.5,10.5) rectangle (11.5,2);

\draw (0,10.5) node[anchor=south,align=left] {\color{green}\Large CPU};
\draw (11,10.5) node[anchor=south,align=right] {\color{red}\Large GPU};

\draw (2.5,9) .. controls (2.5,8.6) and (5.5,9) .. (5.5,8.6);
\draw (8.5,9) .. controls (8.5,8.6) and (5.5,9) .. (5.5,8.6);
\draw[line width=0.7pt] (5.5,8.6) -- (5.5,8.4);
\draw[<-,>=latex] (2.5,8) .. controls (2.5,8.4) and (5.5,8) .. (5.5,8.4) node[midway,above] {$B_{\mu}$};
\draw[<-,>=latex] (8.5,8) .. controls (8.5,8.4) and (5.5,8) .. (5.5,8.4) node[midway,above] {$U$, $V$, $\Sigma$};

\draw (2.5,7) .. controls (2.5,6.6) and (5.5,7) .. (5.5,6.6);
\draw (8.5,7) .. controls (8.5,6.6) and (5.5,7) .. (5.5,6.6);
\draw[line width=0.7pt] (5.5,6.6) -- (5.5,6.4);
\draw[<-,>=latex] (8.5,6) .. controls (8.5,6.4) and (5.5,6) .. (5.5,6.4) node[midway,above] {$W$, $W^{-1}$};

\draw (8.5,4) .. controls (8.5,3.6) and (5.5,4) .. (5.5,3.6);
\draw[line width=0.7pt] (5.5,3.6) -- (5.5,3.4);
\draw[<-,>=latex] (2.5,3) .. controls (2.5,3.4) and (5.5,3) .. (5.5,3.4) node[midway,above] {$A$, $t_{\tau (j)}$};
\end{tikzpicture}
\end{center}
\caption{The hybrid CPU-GPU parallel Prony's algorithm}\label{fig_parallel_prony}
\end{figure}

From the algorithmic point of view, the CPU is dedicated to the more complex tasks of computing SVD, eigendecompositions, and solving the least squares problem, while the GPU is performing matrix-matrix multiplications and summations of matrices. This is so, since multicore CPU, capable of executing several threads simultaneously, excels at computations which require frequent synchronizations, and at the algorithms that require complicated control flows. GPU is optimally exploited for large, massively parallel, and regular computations, where work is equally balanced among GPU cores and where all cores are constantly occupied in full capacity.
With careful choice of the algorithms solving the subtasks, the load between CPU and GPU can be balanced. The basic scheme of the hybrid CPU-GPU algorithm is presented in Figure \ref{fig_parallel_prony}.

The main idea is to minimize communication between CPU and GPU. Therefore, CPU is going to store $T$, $f$, $C_{\mu}$, and $c$, as well as auxiliary matrices $U$, $V$ and $W$, while $B_{\mu}$, $A$, and $t_{\tau (j)}$ are going to be copied from GPU. GPU is going to store $T_{\ell}$, $B_{\mu}$, $S_{\ell}$, $A$, and $t_{\tau (j)}$, while $U$, $V$ and $W$ are going to be copied from CPU. Since $T_{\ell}$ for $\ell =1,\ldots ,d$ are stored on GPU, $S_{\ell}$ have to be generated on GPU as well. On the other side, CPU is diagonalizing $C_{\mu}$, so only this matrix is required, not all $S_{\ell}$. To speed up the computation, we will perform preprocessing step of computing $B_{\mu}=\sum_{\ell =1}^{d}\mu_{\ell}T_{\ell}$ in parallel on GPU which is going to be very fast and simultaneous with SVD computation on CPU. This way only the matrix $B_{\mu}$ will be copied from GPU to CPU, instead of all $S_{\ell}$ matrices.  $C_{\mu}$ is obtained from $B_{\mu}$ in the same way as $S_{\ell}$ from $T_{\ell}$, after $U$, $\Sigma$, and $V$ are obtained.

As we mentioned before, CPU is performing reduced SVD computation of $T$ and diagonalization of $C_{\mu}$, and it is solving the problem of least squares $\text{argmin}_{c}\| A^{T}c-f\|_{2}$. GPU only computes the matrix $B_{\mu}$, performs matrix-matrix multiplications with matrices $T_{\ell}$ and $S_{\ell}$, and generates $z_{\tau (j)}$, $t_{\tau (j)}$ and the matrix $A$.

\subsection{Parallelization on CPU}
Matrix-matrix multiplications are implemented by multithreaded BLAS library \cite{blas}. Only generation of the matrix $T$ is explicitly parallelized, where $T$ is partitioned in uniform block columns, and each CPU thread generates elements in one block column.

\subsection{Parallelization on GPU}
On GPU much more tasks are explicitly parallelized. First of all, each matrix $T_{\ell}$ is built by one block of threads on GPU, where generations of elements is distributed among threads. $B_{\mu}$ is partitioned in uniform block columns, where each block of threads is building one block column, and generation of elements is distributed among threads. Diagonal extraction of each $S_{\ell}$ is performed by one block of threads, where access to the elements of diagonals, $z_{\tau (j)}$ and $t_{\tau (j)}$ is distributed among threads. Generation of $A$ is distributed among threads of different blocks of threads, where each block processes one column. Matrix-matrix multiplications are implemented by cuBLAS \cite{cublas}, a library of parallel BLAS routines implemented on NVIDIA GPUs by CUDA \cite{cuda}. There is a possibility to perform matrix multiplications that generate matrices $S_{\ell}$ in (\ref{eq_generateSl}) and that diagonalize them in (\ref{eq_diagonalizeSl}) in parallel, by exploiting CUDA streams or batched matrix-matrix multiplication kernels. Nevertheless, since the dimensions of involved matrices are rather small the benefit of such approach is negligible.

%----$>$ Provjeriti kako bi funkcioniralo generiranje $S_{\ell}$ it $T_{\ell}$ uz pomo\'{c} \texttt{gemm\_batched()}!

%Treba vidjeti da li se sve osim dijagonalizacije i problema najmanjih kvadrata može riješiti na GPU. Lanczos bi se mogao izvesti na GPU osim onog malog SVD-a, to treba na CPU.

\section{Numerical analysis}\label{sec_num_anal}
One of the important topics of this paper is numerical analysis of the Prony's method in the presence of noisy data. In real applications of the method input data $f(k)$, $k\in I$ may contain some perturbations, for example, due to measurement errors or rounding errors when placing data in computer memory. So, instead of $f(k)$ our input data are
$$\tilde{f}(k)=f(k)+\Delta f_{k},\quad \text{such that  }|\Delta f_{k}|\le \varepsilon |f(k)|,$$
for some relative error $\varepsilon$. Our interest is focused on how size of $\varepsilon$ influences the final result of the algorithm executed in exact arithmetic. Direct consequence of such choice of input data is, that instead of matrices $T$ and $T_{\ell}$, $\ell = 1,\ldots ,d$, we are going to operate with matrices
\begin{equation}\label{eq_tildeT}
\tilde{T}=T+\Delta T,\quad \tilde{T}_{\ell}=T_{\ell}+\Delta T_{\ell},
\end{equation}
such that
$$|\Delta T|\le \varepsilon |T|,\quad |\Delta T_{\ell}|\le \varepsilon |T_{\ell}|,$$
where these inequalities are componentwise, and consequently
\begin{equation}\label{eq_deltaT}
\| \Delta T\|_{F}\le \varepsilon \| T\|_{F},\quad \| \Delta T_{\ell}\|_{F}\le \varepsilon \| T_{\ell}\|_{F},
\end{equation}
where
$\| \cdot \|_{F}$ is Frobenius matrix norm. Since we are computing SVD of $\tilde{T}$, by Hoffman-Wielandt-type bound \cite[Subsection 3.2]{rel_pert_eig_svd} we have bound on errors in singular values
\begin{equation}\label{eq_hoff_wiel}
\sqrt{\sum_{i=1}^{N}|\tilde{\sigma}_{i}-\sigma_{i}|^{2}}\le \| \Delta T\|_{F}\le \varepsilon \| T\|_{F},
\end{equation}
where $\sigma_{i}$ are singular values of $T$, and $\tilde{\sigma}_{i}$ are singular values of $\tilde{T}$, both in nonincreasing order. Let,
\begin{equation}\label{eq_tildeTperp}
\tilde{T}=\left[ \begin{array}{cc} \tilde{U}&\tilde{U}_{\perp} \end{array}\right] \left[ \begin{array}{cc} \tilde{\Sigma}& \\ &\tilde{\Sigma}_{\perp} \end{array}\right] \left[ \begin{array}{c} \tilde{V}^{*}\\ \tilde{V}_{\perp}^{*} \end{array}\right] =\tilde{U}\tilde{\Sigma}\tilde{V}^{*}+\tilde{U}_{\perp}\tilde{\Sigma}_{\perp}\tilde{V}_{\perp}^{*}=\tilde{U}\tilde{\Sigma}\tilde{V}^{*}+T_\perp,
\end{equation}
be SVD of $\tilde{T}$, with $\tilde{U}, \tilde{V}\in \mathbb{C}^{N\times m}$ and $\tilde{\Sigma}\in \mathbb{R}^{m\times m}$ consisting of the $m$ largest singular values $\tilde{\sigma}_{i}$. Further, since the last $N-m$ singular values of $T$ are $0$, from (\ref{eq_hoff_wiel}) it follows
\begin{equation}\label{eq_deltaTperp}
\| T_{\perp}\|_{F}=\sqrt{\sum_{i=m+1}^{N}\tilde{\sigma}_{i}^{2}}\le \varepsilon \| T\|_{F}.
\end{equation}
This bound gives us criteria for numerical rank determination i.e. detection of significant drop in singular values:
\begin{equation}\label{eq_rank_crit}
\tilde{\sigma}_{m+1}\le \sqrt{\sum_{i=m+1}^{N}\tilde{\sigma}_{i}^{2}}\le \varepsilon \| T\|_{F}\le (\varepsilon +\mathcal{O}(\varepsilon^{2}))\| \tilde{T}\|_{F},
\end{equation}
meaning that we can determine numerical rank of $T$ being $m$, when $m$ is the smallest integer satisfying (\ref{eq_rank_crit}). Combining (\ref{eq_tildeT}) and (\ref{eq_tildeTperp}) we obtain
\begin{equation}\label{eq_tildeTsvd}
\tilde{U}\tilde{\Sigma}\tilde{V}^{*}=T+\Delta T-T_\perp = T+\Delta T_{SVD},
\end{equation}
where by (\ref{eq_deltaT}) and (\ref{eq_deltaTperp}) it is
\begin{equation}\label{eq_deltaTsvd}
\| \Delta T_{SVD}\|_{F}\le \| \Delta T\|_{F}+\| T_{\perp}\|_{F}\le 2\varepsilon \| T\|_{F}.
\end{equation}
This is the backward error for SVD, and $\tilde{\tilde{T}}=\tilde{U}\tilde{\Sigma}\tilde{V}^{*}$ is a rank $m$ matrix close to $T$ whose reduced SVD we actually use in the algorithm.

The next step is generation of matrices $S_{\ell}$, $\ell =1,\ldots ,d$. Instead of $S_{\ell}=U^{*}T_{\ell}V\Sigma^{-1}$ we compute
\begin{align}
\tilde{S}_{\ell}&=\tilde{U}^{*}\tilde{T}_{\ell}\tilde{V}\tilde{\Sigma}^{-1}=\tilde{U}^{*}T_{\ell}\tilde{V}\tilde{\Sigma}^{-1}+\tilde{U}^{*}\Delta T_{\ell}\tilde{V}\tilde{\Sigma}^{-1}=\tilde{U}^{*}T_{\ell}\tilde{V}\tilde{\Sigma}^{-1}+\Delta S_{\ell,1}, \label{eq_tildeSl}\\
\| \Delta S_{\ell,1}\|_{F}&\le \frac{1}{\tilde{\sigma}_{m}}\| \Delta T_{\ell}\|_{F}\le \frac{\varepsilon}{\tilde{\sigma}_{m}}\| T_{\ell}\|_{F},\label{eq_deltaSl1}
\end{align}
Further we need to determine a relationship between exact and computed singular vectors. Let $U_{\perp}$ and $V_{\perp}$ be orthonormal matrices whose columns span orthogonal complements of $\text{span}\{ U \}$ and $\text{span}\{ V\}$ respectively. Then we can write
\begin{align}
\tilde{U}&=\left[ \begin{array}{cc} U&U_{\perp} \end{array}\right] \left[ \begin{array}{c} U^{*}\\ U_{\perp}^{*} \end{array}\right] \tilde{U}= U(U^{*}\tilde{U})+U_{\perp}(U_{\perp}^{*}\tilde{U}) \label{eq_tildeU}\\
\tilde{V}&=\left[ \begin{array}{cc} V&V_{\perp} \end{array}\right] \left[ \begin{array}{c} V^{*}\\ V_{\perp}^{*} \end{array}\right] \tilde{V}= V(V^{*}\tilde{V})+V_{\perp}(V_{\perp}^{*}\tilde{V}) \label{eq_tildeV}
\end{align}
Let us assume that all $k$ columns of $X$ are some choice of left singular vectors from $U$, that all $k$ columns of $\tilde{X}$ are some choice of left singular vectors from $\tilde{U}$, and that we have the same assumption for $Y$ and $\tilde{Y}$. By $X_{\perp}$ we denote matrix whose columns are the rest $N-k$ columns of $\left[ \begin{array}{cc} U&U_{\perp} \end{array}\right]$, and we define $Y_{\perp}$ in similar way from $\left[ \begin{array}{cc} V&V_{\perp} \end{array}\right]$. And finally, let $\tau_{i}$ $i=1,\ldots ,k$ are singular values of $T$ that correspond to $X$ and $Y$, and $\tau_{i}$ $i=k+1,\ldots ,N$ are all the rest. $\tilde{\tau}_{i}$ are similarly defined for singular values of $\tilde{\tilde{T}}$. Then, by Wedin's theorem \cite[Theorem 3.3]{rel_pert_eigenspaces}, \cite{wedin_pert_svd}
\begin{equation}
\sqrt{\| X_{\perp}^{*}\tilde{X}\|_{F}^{2}+\| Y_{\perp}^{*}\tilde{Y}\|_{F}^{2}}\le \frac{\sqrt{\| (\tilde{\tilde{T}}-T)Y\|_{F}^{2}+\| (\tilde{\tilde{T}}^{*}-T^{*})X\|_{F}^{2}}}{\delta}\le \frac{\sqrt{2}\| \Delta T_{SVD}\|_{F}}{\delta}\le \frac{2\sqrt{2}\varepsilon \| T\|_{F}}{\delta},
\label{eq_wedintm}
\end{equation}
where
$$\delta =\min \left\{ \min_{i=1,\ldots ,k \atop j=1,\ldots ,N-k}|\tau_{i}-\tilde{\tau}_{k+j}|, \min_{i=1,\ldots ,k}\tau_{i}\right\} .$$
It can be shown that $\| X^{*}\tilde{X}\|_{F}=\| \cos \Theta (X,\tilde{X})\|_{F}$ and $\| X_{\perp}^{*}\tilde{X}\|_{F}=\| \sin \Theta (X,\tilde{X})\|_{F}$, where $\Theta (X,\tilde{X})$ is diagonal matrix with canonical angles between subspaces $\text{span}\{ X\}$ and $\text{span}\{ \tilde{X}\}$. Further, $\cos \Theta (X,\tilde{X})\ge 0$ and $\sin \Theta (X,\tilde{X})\ge 0$.  The same is valid for $Y$ and $\tilde{Y}$.

In case when $X=U(:,i)$, $\tilde{X}=\tilde{U}(:,i)$, $Y=V(:,i)$ and $\tilde{Y}=\tilde{V}(:,i)$ we can conclude the following. Let $\theta_i = \Theta (U(:,i),\tilde{U}(:,i))$, then
$$|(U^{*}\tilde{U})(i,i)|=|U(:,i)^{*}\tilde{U}(:,i)|=\cos \theta_i = \sqrt{1-(\sin \theta_i)^{2}}=1-\frac{1}{2}(\sin \theta_i)^2+\mathcal{O}((\sin \theta_i)^{4}),$$
and by Wedin's theorem  (\ref{eq_wedintm}) we have
$$\sin \theta_i \le \frac{2\sqrt{2}\varepsilon \| T\|_{F}}{\delta_{i}},\quad \delta_{i}=\min \{ \min_{j=1,\ldots ,N \atop j\ne i}|\sigma_{i}-\tilde{\sigma}_{j}|, \sigma_{i}\}.$$
We can scale columns of $U$ by scalars $e^{i\alpha_{i}}$, with absolute value equal to 1, so that $U(:,i)^{*}\tilde{U}(:,i)=\cos \theta_{i}$ for all $i$. The same scaling applies to the columns of $V$. Hence
\begin{align}
(U^{*}\tilde{U})(i,i)&=1+ \Delta U(i,i),\quad |\Delta U(i,i)|\le \frac{4\varepsilon^{2} \| T\|_{F}^{2}}{\delta_{i}^{2}}+\mathcal{O}(\varepsilon^{4}),\label{eq_UUii} \\
(V^{*}\tilde{V})(i,i)&=e^{i\beta_{i}}+ \Delta V_{0}(i,i),\quad |\Delta V_{0}(i,i)|\le \frac{4\varepsilon^{2} \| T\|_{F}^{2}}{\delta_{i}^{2}}+\mathcal{O}(\varepsilon^{4}),\label{eq_VVii}
\end{align}
for some angles $\beta_{i}$. On the other hand,
\begin{equation}\label{eq_UUji}
\sqrt{\sum_{j=1 \atop j\ne i}^{m}|(U^{*}\tilde{U})(j,i)|^{2}}\le \sqrt{\sum_{j=1 \atop j\ne i}^{m}|(U^{*}\tilde{U})(j,i)|^{2}+\sum_{j=1}^{N-m}|(U_{\perp}^{*}\tilde{U})(j,i)|^{2}}=\| X_{\perp}^{*}\tilde{X}\|_{F}\le \frac{2\sqrt{2}\varepsilon \| T\|_{F}}{\delta_{i}}.
\end{equation}
So, from (\ref{eq_UUii}) and (\ref{eq_UUji}) we can conclude that
$$(U^{*}\tilde{U})(:,i)=e_{i}+\Delta U(:,i),\quad \| \Delta U(:,i)\|_{2}^{2}\le \frac{8\varepsilon^{2} \| T\|_{F}^{2}}{\delta_{i}^{2}}+\frac{16\varepsilon^{4} \| T\|_{F}^{4}}{\delta_{i}^{4}}+\mathcal{O}(\varepsilon^{6}).$$
Putting all this together we obtain
\begin{equation}\label{eq_deltaU}
U^{*}\tilde{U}=I+\Delta U,\quad \| \Delta U\|_{F}\le \frac{2\sqrt{2m}\varepsilon \| T\|_{F}}{\delta_{\min}}+\mathcal{O}(\varepsilon^{3}),\ \delta_{\min}=\min_{i=1,\ldots ,m}\delta_{i},
\end{equation}
and in the same way, for $J=\text{diag}(e^{i\beta_{1}},\ldots ,e^{i\beta_{m}})$,
\begin{equation}\label{eq_deltaV0}
V^{*}\tilde{V}=J+\Delta V_{0},\quad \| \Delta V_{0}\|_{F}\le \frac{2\sqrt{2m}\varepsilon \| T\|_{F}}{\delta_{\min}}+\mathcal{O}(\varepsilon^{3}).
\end{equation}

For $X=U$, $\tilde{X}=\tilde{U}$, $Y=V$, and $\tilde{Y}=\tilde{V}$, directly from Wedin's theorem it follows
\begin{align}
\| U_{\perp}^{*}\tilde{U}\|_{F}&\le \frac{2\sqrt{2}\varepsilon \| T\|_{F}}{\delta_{T,\tilde{\tilde{T}}}}, \label{eq_UperpU}\\
\| V_{\perp}^{*}\tilde{V}\|_{F}&\le \frac{2\sqrt{2}\varepsilon \| T\|_{F}}{\delta_{T,\tilde{\tilde{T}}}}, \label{eq_VperpV}\\
\delta_{T,\tilde{\tilde{T}}}&=\min \left\{ \min_{i=1,\ldots ,m \atop j=1,\ldots ,N-m}|\sigma_{i}-\tilde{\sigma}_{m+j}|,\sigma_{m}\right\} \ge \delta_{\min}. \label{eq_delta}
\end{align}

Since from (\ref{eq_hoff_wiel}), (\ref{eq_tildeU}) and (\ref{eq_tildeV}) we have
\begin{align*}
\tilde{\tilde{T}}&=\tilde{U}\tilde{\Sigma}\tilde{V}^{*}=(U(U^{*}\tilde{U})+U_{\perp}(U_{\perp}^{*}\tilde{U}))(\Sigma +\Delta \Sigma)(V(V^{*}\tilde{V})+V_{\perp}(V_{\perp}^{*}\tilde{V}))^{*}\\
&=(U(I+\Delta U)+U_{\perp}(U_{\perp}^{*}\tilde{U}))(\Sigma +\Delta \Sigma)(V(I+(J-I)+\Delta V_{0})+V_{\perp}(V_{\perp}^{*}\tilde{V}))^{*}\\
&=U\Sigma V^{*}+U\Delta U \Sigma V^{*}+U_{\perp}(U_{\perp}^{*}\tilde{U})\Sigma V^{*}+U\Sigma (\bar{J}-I)V^{*}+U\Sigma \Delta V_{0}^{*}V^{*}+\\
&\quad U\Sigma (V_{\perp}^{*}\tilde{V})^{*}V_{\perp}^{*}+ U\Delta \Sigma V^{*}+\mathcal{O}(\varepsilon^{2})\\
&=T+\Delta T_{SVD}
\end{align*}
and $\| \Delta T_{SVD}\|_{F}\le 2\varepsilon \| T\|_{F}$ is small, $\| J-I\|_{F}$ has to be of the same order as $\| \Delta U\|_{F}$, $\| \Delta V_{0}\|_{F}$, $\| U_{\perp}^{*}\tilde{U}\|_{F}$, and $\| V_{\perp}^{*}\tilde{V}\|_{F}$. Hence we can conclude that for $\Delta V=J-I+\Delta V_{0}$
\begin{equation}\label{eq_deltaV}
V^{*}\tilde{V}=I+\Delta V,\quad \| \Delta V\|_{F}\le \frac{g\sqrt{2m}\varepsilon \| T\|_{F}}{\delta_{\min}}+\mathcal{O}(\varepsilon^{3}),
\end{equation}
for some constant $g$.

Now, we switch back to matrices $S_{\ell}$, where by (\ref{eq_tildeSl})
\begin{align*}
S_{\ell}&=U^{*}T_{\ell}V\Sigma^{-1},\quad \ell =1,\ldots, d,\\
\tilde{S}_{\ell}&=\tilde{U}^{*}\tilde{T}_{\ell}\tilde{V}\tilde{\Sigma}^{-1}=\tilde{U}^{*}T_{\ell}\tilde{V}\tilde{\Sigma}^{-1}+\Delta S_{\ell,1},\quad \ell =1,\ldots, d.
\end{align*}
Again, from (\ref{eq_tildeU}) and (\ref{eq_tildeV}) it follows
\begin{align*}
\tilde{S}_{\ell}&=(U(U^{*}\tilde{U})+U_{\perp}(U_{\perp}^{*}\tilde{U}))^{*}T_{\ell}(V(V^{*}\tilde{V})+V_{\perp}(V_{\perp}^{*}\tilde{V}))\tilde{\Sigma}^{-1}+\Delta S_{\ell,1}\\
&=(\tilde{U}^{*}U)[U^{*}T_{\ell}V\Sigma^{-1}]\Sigma (V^{*}\tilde{V})\tilde{\Sigma}^{-1}+(\tilde{U}^{*}U_{\perp})U_{\perp}^{*}T_{\ell}V(V^{*}\tilde{V})\tilde{\Sigma}^{-1}+\\
&\quad +(\tilde{U}^{*}U)U^{*}T_{\ell}V_{\perp}(V_{\perp}^{*}\tilde{V})\tilde{\Sigma}^{-1}+(\tilde{U}^{*}U_{\perp})U_{\perp}^{*}T_{\ell}V_{\perp}(V_{\perp}^{*}\tilde{V})\tilde{\Sigma}^{-1} +\Delta S_{\ell,1}
\end{align*}
If we define
\begin{align*}
\Delta S_{\ell,2}&=(\tilde{U}^{*}U_{\perp})U_{\perp}^{*}T_{\ell}V(V^{*}\tilde{V})\tilde{\Sigma}^{-1},\\
\Delta S_{\ell,3}&=(\tilde{U}^{*}U)U^{*}T_{\ell}V_{\perp}(V_{\perp}^{*}\tilde{V})\tilde{\Sigma}^{-1},\\
\Delta S_{\ell,4}&=(\tilde{U}^{*}U_{\perp})U_{\perp}^{*}T_{\ell}V_{\perp}(V_{\perp}^{*}\tilde{V})\tilde{\Sigma}^{-1},
\end{align*}
and if we take into account (\ref{eq_deltaU}) and (\ref{eq_deltaV}) we obtain
\begin{align*}
\tilde{S}_{\ell}&=(I+\Delta U^{*})S_{\ell}\Sigma (I+\Delta V)\tilde{\Sigma}^{-1}+\Delta S_{\ell,1}+\Delta S_{\ell,2}+\Delta S_{\ell,3}+\Delta S_{\ell,4}\\
&=S_{\ell}+S_{\ell}(\Sigma \tilde{\Sigma}^{-1}-I)+\Delta U^{*}S_{\ell}\Sigma \tilde{\Sigma}^{-1}+S_{\ell}\Sigma \Delta V\tilde{\Sigma}^{-1}+\Delta U^{*}S_{\ell}\Sigma \Delta V\tilde{\Sigma}^{-1}+\\
&\quad +\Delta S_{\ell,1}+\Delta S_{\ell,2}+\Delta S_{\ell,3}+\Delta S_{\ell,4}.
\end{align*}
Further, if we define
\begin{align*}
\Delta S_{\ell,5}&=\Delta U^{*}S_{\ell}\Sigma \tilde{\Sigma}^{-1}+S_{\ell}\Sigma \Delta V\tilde{\Sigma}^{-1}+\Delta U^{*}S_{\ell}\Sigma \Delta V\tilde{\Sigma}^{-1},\\
\Delta S_{\ell,6}&=S_{\ell}(\Sigma \tilde{\Sigma}^{-1}-I),\\
\Delta S_{\ell}&=\Delta S_{\ell,1}+\Delta S_{\ell,2}+\Delta S_{\ell,3}+\Delta S_{\ell,4}+\Delta S_{\ell,5}+\Delta S_{\ell,6}
\end{align*}
it only remains to find norm bounds on $\Delta S_{\ell,k}$ for $k=2,\ldots ,6$.
\begin{align*}
\| \Delta S_{\ell,2}\|_{2}&\le \| \tilde{U}^{*}U_{\perp}\|_{F}\| T_{\ell}\|_{F}\| I+\Delta V\|_{2}\| \tilde{\Sigma}^{-1}\|_{2}\le \frac{2\sqrt{2}\varepsilon \| T\|_{F}\| T_{\ell}\|_{F}}{\delta_{T,\tilde{\tilde{T}}}\tilde{\sigma}_{m}}\left( 1+\frac{g\sqrt{2m}\varepsilon \| T\|_{F}}{\delta_{\min}}+\mathcal{O}(\varepsilon^{3})\right) \\
&\le \frac{2\sqrt{2}\varepsilon \| T\|_{F}\| T_{\ell}\|_{F}}{\delta_{min}\tilde{\sigma}_{m}}+\mathcal{O}(\varepsilon^{2}),\\
\| \Delta S_{\ell,3}\|_{F}&\le \| I+\Delta U \|_{2}\| T_{\ell}\|_{F}\| V_{\perp}^{*}\tilde{V}\|_{F}\| \tilde{\Sigma}^{-1}\|_{2}\le \frac{2\sqrt{2}\varepsilon \| T\|_{F}\| T_{\ell}\|_{F}}{\delta_{min}\tilde{\sigma}_{m}}+\mathcal{O}(\varepsilon^{2}),\\
\| \Delta S_{\ell,4}\|_{F}&\le \| \tilde{U}^{*}U_{\perp}\|_{F}\| T_{\ell}\|_{F}\| V_{\perp}^{*}\tilde{V}\|_{F}\| \tilde{\Sigma}^{-1}\|_{2}\le \mathcal{O}(\varepsilon^{2}),\\
\| \Delta S_{\ell,5}\|_{F}&\le \| \Delta U\|_{F}\| S_{\ell}\Sigma\|_{F} \| \tilde{\Sigma}^{-1}\|_{2}+\| S_{\ell}\Sigma\|_{F} \| \Delta V\|_{F}\| \tilde{\Sigma}^{-1}\|_{2}+\| \Delta U^{*}\|_{F}\| S_{\ell}\Sigma\|_{F} \| \Delta V\|_{F}\| \tilde{\Sigma}^{-1}\|_{2}\\
&\le \frac{(2+g)\sqrt{2m}\varepsilon \| T\|_{F}\| T_{\ell}\|_{F}}{\delta_{\min}\tilde{\sigma}_{m}}+\mathcal{O}(\varepsilon^{2}),\\
\| \Delta S_{\ell,6}\|_{F}&\le \| S_{\ell}\|_{F}\left\| \text{diag}\left( \frac{\sigma_{i}-\tilde{\sigma}_{i}}{\tilde{\sigma}_{i}}\right) \right\|_{F}\le \frac{\varepsilon \| T\|_{F}}{\tilde{\sigma}_{m}}\| S_{\ell}\|_{F}\le \frac{\varepsilon \| T\|_{F}\| T_{\ell}\|_{F}}{\sigma_{m}\tilde{\sigma}_{m}}.
\end{align*}
Finally we can conclude that
\begin{align}
\tilde{S}_{\ell}&=S_{\ell}+\Delta S_{\ell}, \nonumber\\
\| \Delta S_{\ell}\|_{F}&\le \frac{\varepsilon \| T_{\ell}\|_{F}}{\tilde{\sigma}_{m}}\left( 1+ \frac{4\sqrt{2}\| T\|_{F}}{\delta_{min}}+ \frac{(2+g)\sqrt{2m}\| T\|_{F}}{\delta_{\min}}+ \frac{\| T\|_{F}}{\sigma_{m}}\right) + \mathcal{O}(\varepsilon^{2}) \nonumber\\
&\le \frac{\varepsilon \| T_{\ell}\|_{F}}{\tilde{\sigma}_{m}}\left( 1+\frac{(1+4\sqrt{2}+(2+g)\sqrt{2m})\| T \|_{F}}{\delta_{\min}}\right) +\mathcal{O}(\varepsilon^{2}). \label{eq_deltaSl}
\end{align}

Now we switch to $C_{\mu}$, and we will operate with
\begin{equation}\label{eq_Cmu}
\tilde{C}_{\mu}=\sum_{\ell =1}^{d}\mu_{\ell}\tilde{S}_{\ell}=\sum_{\ell =1}^{d}\mu_{\ell}S_{\ell}+\sum_{\ell =1}^{d}\mu_{\ell}\Delta S_{\ell}=C_{\mu}+\Delta C_{\mu}, \quad \| \mu \|_{2}=1
\end{equation}
instead, where
\begin{align}
\| \Delta C_{\mu}\|_{F}&\le \max_{\ell =1,\ldots ,d}\| \Delta S_{\ell}\|_{F}\sum_{\ell =1}^{d}\mu_{\ell} \nonumber \\
&\le \frac{\sqrt{d}\varepsilon \max_{\ell =1,\ldots ,d}\| T_{\ell}\|_{F}}{\tilde{\sigma}_{m}}\left( 1+\frac{(1+4\sqrt{2}+(2+g)\sqrt{2m})\| T \|_{F}}{\delta_{\min}}\right) +\mathcal{O}(\varepsilon^{2}).\label{eq_deltaCmu}
\end{align}
We will turn now to the perturbation theory for eigenvalue problem, and we will investigate relation between spectral decompositions of $C_{\mu}$ and $\tilde{C}_{\mu}$. Let us define spectral decompositions
\begin{align}
C_{\mu}&=W\Lambda_{\mu}W^{-1}, &\Lambda_{\mu}= \text{diag}(\lambda_{1},\ldots ,\lambda_{m}) \label{eq_Cmu_diag}\\
\tilde{C}_{\mu}&=\tilde{W}\tilde{\Lambda}_{\mu}\tilde{W}^{-1}, &\tilde{\Lambda}_{\mu}= \text{diag}(\tilde{\lambda}_{1},\ldots ,\tilde{\lambda}_{m}) \label{eq_tildeCmu_diag}
\end{align}
If we assume that all eigenvalues of $C_{\mu}$ and $\tilde{C}_{\mu}$ are distinct and ordered in ascending orther, then columns of $W=[\begin{array}{ccc} w_{1},\ldots ,w_{m} \end{array}]$ and $\tilde{W}=[\begin{array}{ccc} \tilde{w}_{1},\ldots ,\tilde{w}_{m} \end{array}]$ are corresponding eigenvectors which are usually normalized to have norm 1. Further, we choose again to scale the vectors $w_{1}$,\ldots, $w_{m}$ by $e^{i\alpha_{i}}$ so that $\tilde{w}_{i}^{*}w_{i}=\cos \phi_{i}$ and not only $|\tilde{w}_{i}^{*}w_{i}|=\cos \phi_{i}$, where $\phi_{i}$ is the angle between eigenvectors $w_{i}$ and $\tilde{w}_{i}$ belonging to the eigenvalues $\lambda_{i}$ and $\tilde{\lambda}_{i}$. Then,
$$\| \tilde{w}_{i}-w_{i}\|_{2}=\sqrt{\| \tilde{w}_{i}\|_{2}^{2}+\| w_{i}\|_{2}^{2}-2\cos \phi_{i} \| \tilde{w}_{i}\|_{2}\| w_{i}\|_{2}}=\sqrt{2}\sqrt{1-\cos \phi_{i}}.$$
By \cite[Theorem 5.1 ]{abs_rel_pert_inv_subsp} we can conclude that
$$|\sin \phi_{i}|\le \frac{\kappa_{2}(W^{-1}(j\ne i,:))\kappa_{2}(\tilde{w}_{i})\| \Delta C_{\mu}\|_{F}}{\min_{j=1,\ldots ,m \atop j\ne i}|\lambda_{j}-\tilde{\lambda}_{i}|}.$$
By interlacing property for singular values \cite[Corollary 3.1.3]{horn_johnson_topics} it follows that
$$\kappa_{2}(W^{-1}(j\ne i,:))\le \kappa_{2}(W^{-1})=\kappa_{2}(W),$$
and $\kappa_{2}(\tilde{w}_{i})=1$, so finally
\begin{equation}\label{eq_sin_phi_i}
|\sin \phi_{i}|\le \frac{\kappa_{2}(W)\| \Delta C_{\mu}\|_{F}}{\min_{j=1,\ldots ,m \atop j\ne i}|\lambda_{j}-\tilde{\lambda}_{i}|}.
\end{equation}
Now we can finish the bound on $\tilde{w}_{i}-w_{i}$
\begin{equation}\label{eq_deltaw_i}
\| \tilde{w}_{i}-w_{i}\|_{2}=\sqrt{2}\sqrt{1-\sqrt{1-(\sin \phi_{i})^2}}=\sqrt{\frac{(\sin \phi_{i})^2}{1+\sqrt{1-(\sin \phi_{i})^2}}}\le |\sin \phi_{i}|.
\end{equation}
Thus
\begin{equation}\label{eq_tildeW}
\tilde{W}=W+\Delta W
\end{equation}
where
\begin{align}
\| \Delta W\|_{F}&=\sqrt{\sum_{i=1}^{m}\| \tilde{w}_{i}-w_{i}\|_{2}^{2}} \nonumber \\
&\le \frac{\sqrt{m}\kappa_{2}(W)}{\min_{i,j=1,\ldots ,m \atop j\ne i}|\lambda_{j}-\tilde{\lambda}_{i}|}\cdot
\frac{\sqrt{d}\varepsilon \max_{\ell =1,\ldots ,d}\| T_{\ell}\|_{F}}{\tilde{\sigma}_{m}}\left( 1+\frac{(1+4\sqrt{2}+(2+g)\sqrt{2m})\| T \|_{F}}{\delta_{\min}}\right) +\mathcal{O}(\varepsilon^{2}) \nonumber \\
&=\frac{\sqrt{md}\varepsilon \max_{\ell =1,\ldots ,d}\| T_{\ell}\|_{F}\kappa_{2}(W)}{\gamma \tilde{\sigma}_{m}}\left( 1+\frac{(1+4\sqrt{2}+(2+g)\sqrt{2m})\| T \|_{F}}{\delta_{\min}}\right) +\mathcal{O}(\varepsilon^{2}) \label{eq_deltaW}
\end{align}
where
\begin{equation}\label{eq_gamma}
\gamma=\min_{i,j=1,\ldots ,m \atop j\ne i}|\lambda_{j}-\tilde{\lambda}_{i}|.
\end{equation}
We also need a bound on difference $\tilde{W}^{-1}-W^{-1}$. Since,
$$\tilde{W}=W(I+W^{-1}\Delta W),$$
in case when
$$\| W^{-1}\Delta W\|_{F}\le \frac{\| \Delta W\|_{F}}{\sigma_{min}(W)}<1,$$
where $\sigma_{min}(W)$ is the smallest singular value of $W$, we have
\begin{equation}\label{eq_tildeWinv}
\tilde{W}^{-1}=\left( I-W^{-1}\Delta W+(W^{-1}\Delta W)^2-(W^{-1}\Delta W)^3+\cdots \right) W^{-1}.
\end{equation}

If we define by $S_{\ell}=W\Lambda_{\ell}W^{-1}$ spectral decomposition of $S_{\ell}$, then we further compute
\begin{align*}
\tilde{W}^{-1}\tilde{S}_{\ell}\tilde{W}&=\left( I-W^{-1}\Delta W+(W^{-1}\Delta W)^2-\cdots \right) W^{-1}(S_{\ell}+\Delta S_{\ell})W(I+W^{-1}\Delta W)\\
&=\left( I-W^{-1}\Delta W+(W^{-1}\Delta W)^2-\cdots \right)\Lambda_{\ell}(I+W^{-1}\Delta W)+W^{-1}\Delta S_{\ell}W+\mathcal{O}(\varepsilon^{2})\\
&=\Lambda_{\ell}-W^{-1}\Delta W\Lambda_{\ell}+ \Lambda_{\ell}W^{-1}\Delta W+W^{-1}\Delta S_{\ell}W+\mathcal{O}(\varepsilon^{2})\\
&=\Lambda_{\ell}+\Delta \Lambda_{\ell}
\end{align*}
where
\begin{align}
\| \Delta \Lambda_{\ell}\|_{F}&\le 2\frac{\| \Delta W\|_{F}}{\sigma_{min}(W)}\| \Lambda_{\ell}\|_{2}+\kappa_{2}(W)\| \Delta S_{\ell}\|_{F}+\mathcal{O}(\varepsilon^{2}) \nonumber\\
&\le \varepsilon \left( \frac{2\sqrt{md} \| \Lambda_{\ell}\|_{2}}{\gamma \tilde{\sigma}_{m}\sigma_{min}(W)}+ \frac{1}{\tilde{\sigma}_{m}}\right)\left( 1+\frac{(1+4\sqrt{2}+(2+g)\sqrt{2m})\| T \|_{F}}{\delta_{\min}}\right) \max_{\ell =1,\ldots ,d}\| T_{\ell}\|_{F}\kappa_{2}(W) +\mathcal{O}(\varepsilon^{2}).\label{eq_deltaLambda_l}
\end{align}
Since we are extracting the diagonal of $\tilde{W}^{-1}\tilde{S}_{\ell}\tilde{W}$, we can bound it from the exact eigenvalues of $S_{\ell}$ by
\begin{equation}\label{eq_diag_tildeSl}
\| \text{diag}(\tilde{W}^{-1}\tilde{S}_{\ell}\tilde{W}-\Lambda_{\ell})\|_{2}\le \| \Delta \Lambda_{\ell}\|_{F}
\end{equation}
where $\| \Delta \Lambda_{\ell}\|_{F}$ is bounded by (\ref{eq_deltaLambda_l}), and $\Lambda_{\ell}=\text{diag}(z_{1}(\ell),\ldots,z_{m}(\ell))$ with $z_{j}(\ell)$ denoting $\ell$-th component of $z_{j}$. Furthermore,
\begin{equation}\label{eq_tildezj}
\left( \tilde{W}^{-1}\tilde{S}_{\ell}\tilde{W}\right) (j,j)=\tilde{z}_{j}(\ell),\quad \tilde{z}_{j}=z_{j}+\Delta z_{j},
\end{equation}
and each component of $\Delta z_{j}$ is bounded by (\ref{eq_deltaLambda_l}). Further, $\tilde{t}_{j}$ are computed from $\tilde{z}_{j}$ as
\begin{equation}\label{eq_tildetj}
\tilde{t}_{j}(i)=\text{Re}\left( \frac{\log \tilde{z}_{j}(i)}{-2\pi \iota}\right) =\frac{\text{Im}(\log \tilde{z}_{j}(i))}{2\pi}=\frac{1}{2\pi}\text{arctg} \frac{\text{Im}(\tilde{z}_{j}(i))}{\text{Re}(\tilde{z}_{j}(i))} .
\end{equation}
If we analyze this expression we obtain
$$
\text{arctg}\frac{\text{Im}(\tilde{z}_{j}(i))}{\text{Re}(\tilde{z}_{j}(i))} = \text{arctg}\frac{\text{Im}(z_{j}(i))+\text{Im}(\Delta z_{j}(i))}{\text{Re}(z_{j}(i))+\text{Re}(\Delta z_{j}(i))}
 =\text{arctg}\frac{\text{Im}(z_{j}(i))}{\text{Re}(z_{j}(i))}\left( \frac{1+\frac{\text{Im}(\Delta z_{j}(i))}{\text{Im}(z_{j}(i))}}{1+\frac{\text{Re}(\Delta z_{j}(i))}{\text{Re}(z_{j}(i))}}\right)
$$
If we denote by
$$q_{j}(i)=\frac{\text{Im}(z_{j}(i))}{\text{Re}(z_{j}(i))},\quad \tilde{q}_{j}(i)=\frac{\text{Im}(\tilde{z}_{j}(i))}{\text{Re}(\tilde{z}_{j}(i))},\quad \delta p_{j}(i)= \frac{\text{Im}(\Delta z_{j}(i))}{\text{Im}(z_{j}(i))},\quad \delta r_{j}(i)=\frac{\text{Re}(\Delta z_{j}(i))}{\text{Re}(z_{j}(i))},$$
then we have
$$
\text{arctg }\tilde{q}_{j}(i)=\text{arctg }\left( q_{j}(i)\frac{1+\delta p_{j}(i)}{1+\delta r_{j}(i)}\right) ,
$$
with
\begin{align}
\frac{1+\delta p_{j}(i)}{1+\delta r_{j}(i)}&=1+\frac{\delta p_{j}(i)-\delta r_{j}(i)}{1+\delta r_{j}(i)}=1+\delta q_{j}(i), \nonumber \\
|\delta q_{j}(i)|&\le \frac{2\sqrt{\delta r_{j}(i)^2+\delta p_{j}(i)^2}}{1-\sqrt{\delta r_{j}(i)^2+\delta p_{j}(i)^2}},\quad \text{where} \nonumber\\
\sqrt{\delta r_{j}(i)^2+\delta p_{j}(i)^2}&=\sqrt{\frac{\text{Re}(\Delta z_{j}(i))^{2}\text{Im}(z_{j}(i))^{2}+\text{Im}(\Delta z_{j}(i))^{2}\text{Re}(z_{j}(i))^{2}}{\text{Re}(z_{j}(i))^{2}\text{Im}(z_{j}(i))^{2}}}\le \frac{1}{\eta} |\Delta z_{j}(i)| \nonumber
\end{align}
%\frac{2|\Delta z_{j}(i)|}{1-|\Delta z_{j}(i)|}=2|\Delta z_{j}(i)|+\mathcal{O}(\varepsilon^{2}),\label{eq_deltaqj}
for $\eta =\min_{j=1,\ldots ,m \atop i=1,\ldots ,d}|\text{Re}(z_{j}(i))\text{Im}(z_{j}(i))|$ and taking into account $|z_{j}(i)|=1$. Thus
\begin{equation}\label{eq_deltaqj}
|\delta q_{j}(i)|\le \frac{2|\Delta z_{j}(i)|}{\eta-|\Delta z_{j}(i)|}=\frac{2|\Delta z_{j}(i)|}{\eta}+\mathcal{O}(\varepsilon^{2})
\end{equation}
If we assume that for $|q_{j}(i)|<1$ it is $|\tilde{q}_{j}(i)|<1$, and $|\tilde{q}_{j}(i)|>1$ for $|q_{j}(i)|>1$, from Taylor expansion of function arctg we can conclude that
\begin{equation}\label{eq_deltaxij}
\text{arctg }\tilde{q}_{j}(i)=\text{arctg }q_{j}(i)+\xi_{j}(i), \quad |\xi_{j}(i)|\le |\delta q_{j}(i)|+\mathcal{O}(\varepsilon^{2}),
\end{equation}
and from (\ref{eq_deltaqj})
\begin{equation}\label{eq_tildetj}
\tilde{t}_{j}(i)=t_{j}(i)+\Delta t_{j}(i),\quad |\Delta t_{j}(i)|\le \frac{|\Delta z_{j}(i)|}{\pi \eta }+\mathcal{O}(\varepsilon^{2}),
\end{equation}
where $|\Delta z_{j}(i)|$ is bounded by (\ref{eq_deltaLambda_l}).

The only thing that remains is to analyze solution of the least squares problem
\begin{equation}\label{eq_tilde_LS}
\min_{\tilde{c}}\| \tilde{A}^{T}\tilde{c}-\tilde{f}\|_{2},
\end{equation}
where $\tilde{A}=[\tilde{z}_{j}^{k}]_{j=1,\ldots,m,\ k\in I_{n}}\in \mathbb{C}^{m\times N}$. Then we have
\begin{align}
\tilde{z}_{j}^{k}&=\tilde{z}_{j}(1)^{k(1)}\cdots \tilde{z}_{j}(d)^{k(d)}=(z_{j}(1)+\Delta z_{j}(1))^{k(1)}\cdots (z_{j}(d)+\Delta z_{j}(d))^{k(d)} \nonumber \\
&=\left( z_{j}(1)^{k(1)}+k(1)z_{j}(1)^{k(1)-1}\Delta z_{j}(1)+\mathcal{O}(\varepsilon^{2})\right) \cdots \left( z_{j}(d)^{k(d)}+k(d)z_{j}(d)^{k(d)-1}\Delta z_{j}(1)+\mathcal{O}(\varepsilon^{2})\right) \nonumber \\
&=z_{j}(1)^{k(1)}\cdots z_{j}(d)^{k(d)}+z_{j}(1)^{k(1)}\cdots z_{j}(d)^{k(d)}\left( \frac{k(1)\Delta z_{j}(1)}{z_{j}(1)}+\cdots +\frac{k(d)\Delta z_{j}(d)}{z_{j}(d)}\right) +\mathcal{O}(\varepsilon^{2}) \nonumber \\
&=z_{j}^{k}(1+\Delta z_{j}^{k}) \label{eq_tildezjk}
\end{align}
with
\begin{equation}\label{eq_deltazjk}
|\Delta z_{j}^{k}|\le nd\max_{\ell =1,\ldots ,d}\| \Delta \Lambda_{\ell}\|_{F}+\mathcal{O}(\varepsilon^{2}).
\end{equation}
Hence, the following componentwise bound holds
\begin{align}
\tilde{A}&=A+\Delta A, \label{eq_tildeA}\\
|\Delta A|&\le nd\max_{\ell =1,\ldots ,d}\| \Delta \Lambda_{\ell}\|_{F}|A|+\mathcal{O}(\varepsilon^{2}) \label{eq_deltaA}
\end{align}
and
\begin{align}
\| \Delta A\|_{F}&=\sqrt{\sum_{j=1}^{m}\sum_{k\in I_{n}}|\Delta z_{j}^{k}|^{2}|z_{j}^{k}|^{2}}\le \sqrt{mN}nd\max_{\ell =1,\ldots ,d}\| \Delta \Lambda_{\ell}\|_{F}\| A\|_{F}+\mathcal{O}(\varepsilon^{2})\nonumber \\
&\le \varepsilon \sqrt{mN}nd\left( \frac{2\sqrt{md} \| \Lambda_{\ell}\|_{2}}{\gamma \tilde{\sigma}_{m}\sigma_{min}(W)}+ \frac{1}{\tilde{\sigma}_{m}}\right)\left( 1+\frac{(1+4\sqrt{2}+(2+g)\sqrt{2m})\| T \|_{F}}{\delta_{\min}}\right) \cdot \nonumber\\
&\quad \cdot \max_{\ell =1,\ldots ,d}\| T_{\ell}\|_{F}\kappa_{2}(W)\| A\|_{F} +\mathcal{O}(\varepsilon^{2}) \label{eq_deltaA_F}
\end{align}
or simpler
\begin{align}
\| \Delta A\|_{F}&=\zeta \| A\|_{F}, \nonumber\\
\zeta&\le \mathcal{O}\left( \frac{d^{\frac{3}{2}}m^{\frac{3}{2}}nN^{\frac{1}{2}}}{\gamma \delta_{\min} \tilde{\sigma}_{m}\sigma_{min}(W)}\right) \| T \|_{F} \max_{\ell =1,\ldots ,d}\| T_{\ell}\|_{F}\kappa_{2}(W)\varepsilon +\mathcal{O}(\varepsilon^{2}).\label{eq_eta}
\end{align}

Finally, we employ perturbation theory for linear least squares problems, where for the final result we need the following
\begin{align*}
\tilde{f}&=f+\Delta f, &\| \Delta f\|_{2}\le \varepsilon \| f\|_{2}\\
\tilde{A}&=A+\Delta A, &\| \Delta A\|_{2}\le \| \Delta A\|_{F}=\zeta \| A\|_{F}\le \sqrt{m}\zeta \| A\|_{2}\\
r&=f-A^{T}c=0
\end{align*}
then, by \cite[Theorem 20.1]{accur_stab} provided that $\kappa_{2}(A) \sqrt{m}\zeta < 1$ it follows
\begin{equation}\label{eq_deltac}
\frac{\| \tilde{c}-c\|_{2}}{\| c\|_{2}}\le \frac{2\kappa_{2}(A)\sqrt{m}\zeta}{1-\kappa_{2}(A)\sqrt{m}\zeta}\le \mathcal{O}\left( \frac{d^{\frac{3}{2}}m^{2}nN^{\frac{1}{2}}}{\gamma \delta_{\min} \tilde{\sigma}_{m}\sigma_{min}(W)}\right) \| T \|_{F} \max_{\ell =1,\ldots ,d}\| T_{\ell}\|_{F}\kappa_{2}(W)\kappa_{2}(A)\varepsilon +\mathcal{O}(\varepsilon^{2}).
\end{equation}

As we can see from derived analysis the forward errors of computed parameters $\tilde{t}_{j}$ and coefficients $\tilde{c}_{j}$ are proportional to the relative error $\varepsilon$ of the input data, and depend on condition numbers of the eigenvalue matrix $W$ and least square problem matrix $A$, as expected from the perturbation theory. The gap function on singular values of $T$ and $\tilde{T}$ also plays an important role in the error bound. In case when the algorithm is executed in finite precision arithmetics there would be additional terms in error bounds proportional to the unit roundoff.

%prokomentirati konačne rezultate o čemu ovise greške za t i c.

\section{Numerical tests}
The aim of this section is to prove efficiency of the parallel variant of Prony's method, as well as to illustrate results of its numerical analysis.

In the efficiency tests the sampling values were generated from a predetermined exponential sum with no noise, where $m$ was chosen from the set $\{ 5,10,15,20\}$, $d$ was set to $2$ and $3$. The sample size was fixed for $n=20$, which produces matrix dimension of $N=441$ for $d=2$, and $N=9261$ for $d=3$. The parameters and coefficients are defined as follows:
$$
t_{j}(i)=((i-1)\cdot m+j-1)\cdot 10^{-\lceil \log_{10}(d\cdot m)\rceil},\quad c_{j}=j+\iota j,\quad i=1,\ldots ,d,\ j=1,\ldots ,m.
$$
The criteria for the rank determination was the first $i$ that satisfies $\sigma_{i}<N\varepsilon_{M}\sigma_{1}$, where $\varepsilon_{M}$ is the machine precision. This tolerance is chosen in order to incorporate possible rounding error of the SVD algorithm. The following computational environment was used:
\begin{itemize}
	\item 2x Intel(R) Xeon(R) E5-2690 v3 @ 2.60GHz (24 cores in total);
	\item 256 GB RAM, each processor is equipped with 30 MB of cache memory;
	\item Nvidia Tesla K40c (Kepler generation, 12 GB of GDDR5);
	\item Intel Parallel Studio XE 2016 + MKL 11.3;
	\item Nvidia CUDA 8.0.
    \item The CPUs reach the peak DGEMM performance of about $800$ Gflops, while the GPU has the peak performance of about $1200$ Gflops.
\end{itemize}

\begin{figure}[htbp]
  \centering
  \includegraphics[width=0.47\textwidth, trim=100 245 100 240, clip]{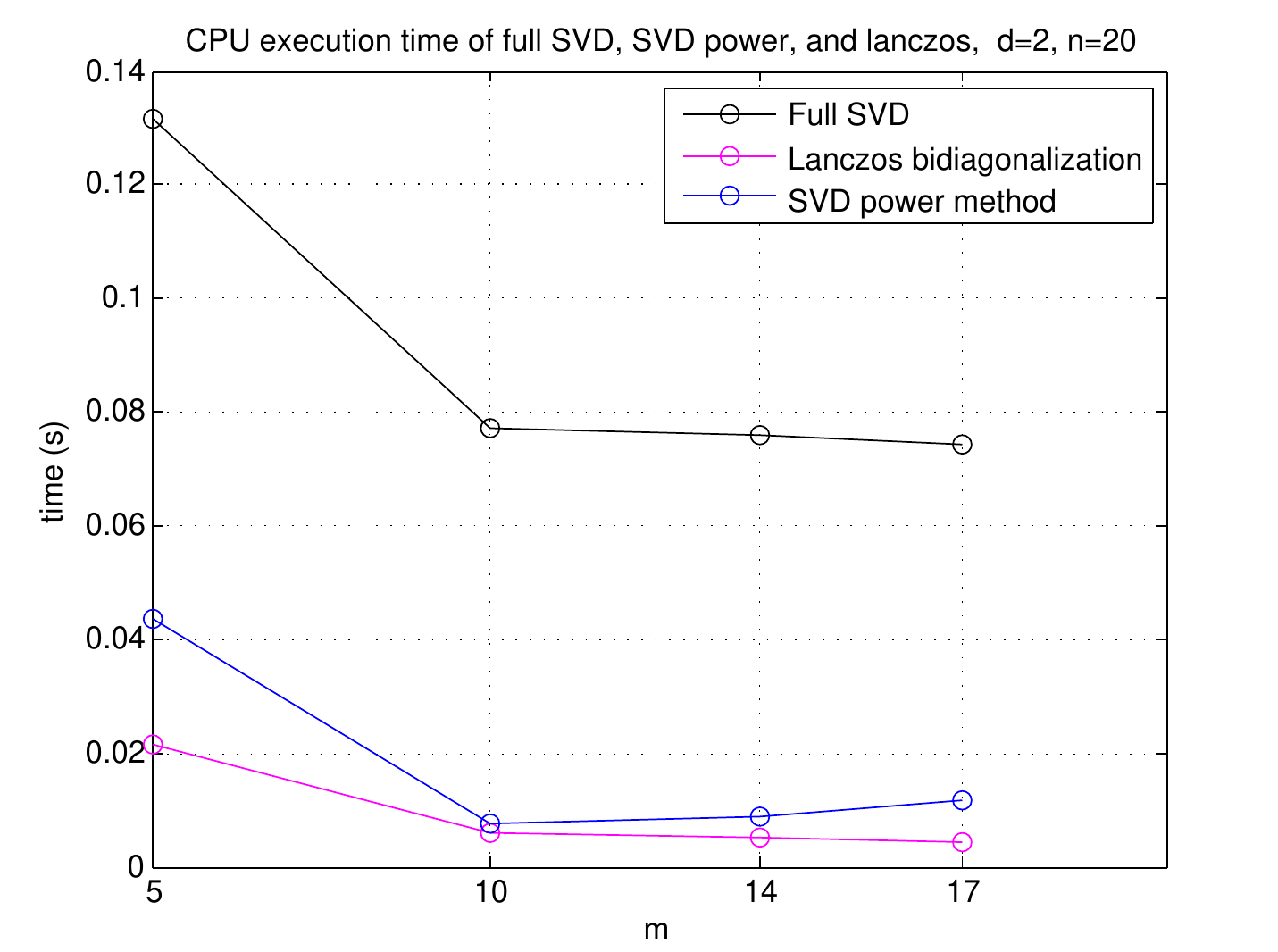} \includegraphics[width=0.47\textwidth, trim=100 245 100 240, clip]{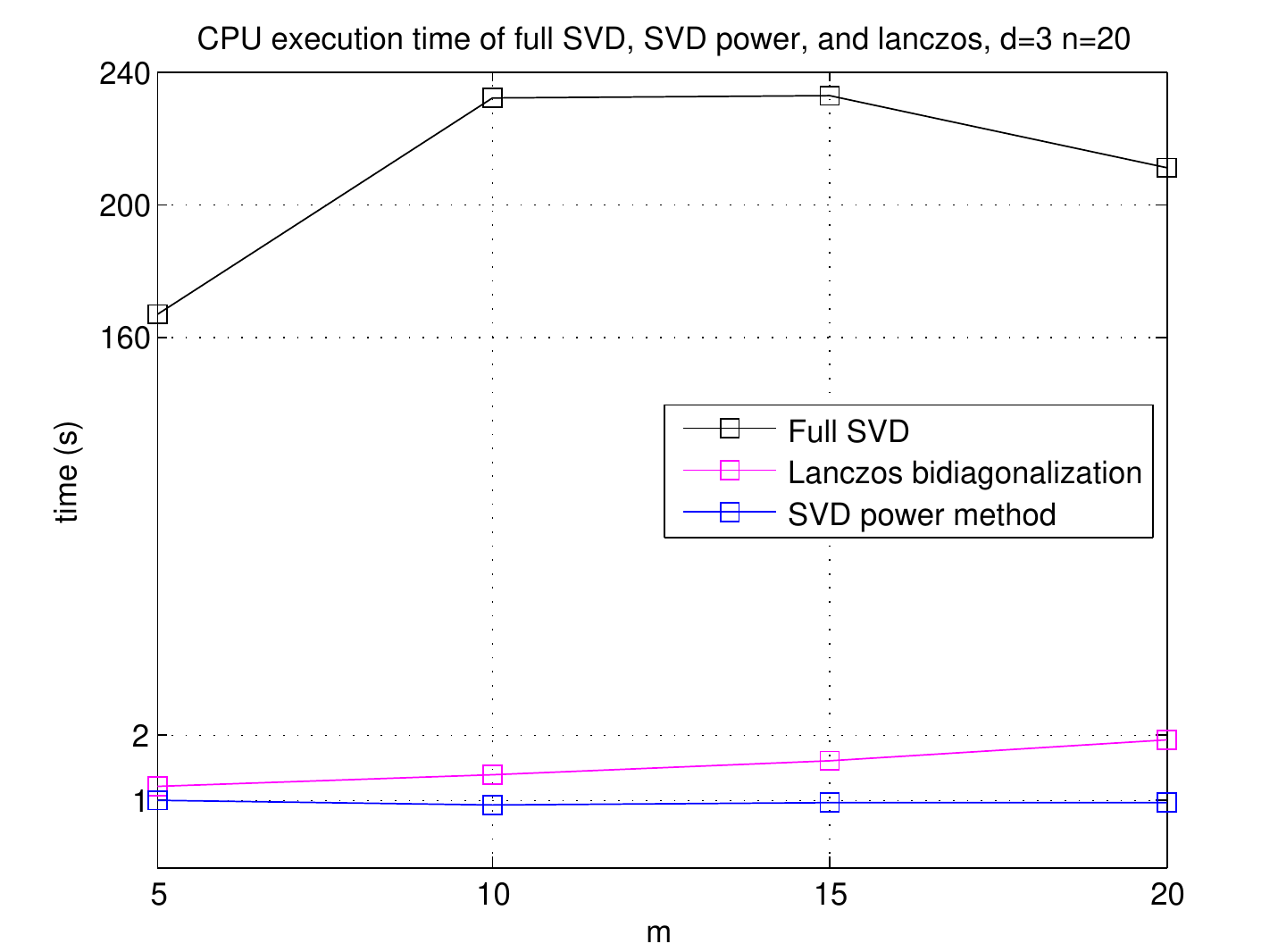}
  \caption{CPU execution time of three SVD algorithms, for $d=2$ and $d=3$.}\label{fig_time_svd}
\end{figure}

First, we are going to compare the SVD algorithms. SVD computation of the matrix $T$ is performed on CPU for all algorithm variants, hence we measured CPU execution time for three SVD algorithms:
\begin{enumerate}
\item computation of full SVD implemented by LAPACK function \texttt{zgesvd()}
\item Lanczos bidiagonalization method
\item block power method for SVD, with starting column-dimensions of matrices $U_{0}$ and $V_{0}$ equal to $2m$.
\end{enumerate}
The execution times are presented in Figure \ref{fig_time_svd}.
It is interesting to see that all three algorithms determined the same rank, and in cases when $d=2$ and $m=15,20$ it seems that condition on $n$ from \cite[Theorem 2.1]{prony_matrix} was not satisfied and sample was too small. Determined rank of $T$ was smaller than the number of terms in the exponential sum of the sample. In these cases the accuracy of the Prony's solution was reduced. Further, we can also conclude that for smaller matrix dimensions Lanczos bidiagonalization method is the fastest method, but for larger dimensions the SVD power method is the most efficient as expected, specially for larger ranks. Computation of full SVD is by far the slowest method, even 200 times slower than the power method for $N=9261$. Thus, we can assert that right choice of the SVD algorithm is crucial for efficiency of the Prony's method.

\begin{figure}[htbp]
  \centering
  \includegraphics[width=0.47\textwidth, trim=100 245 100 240, clip]{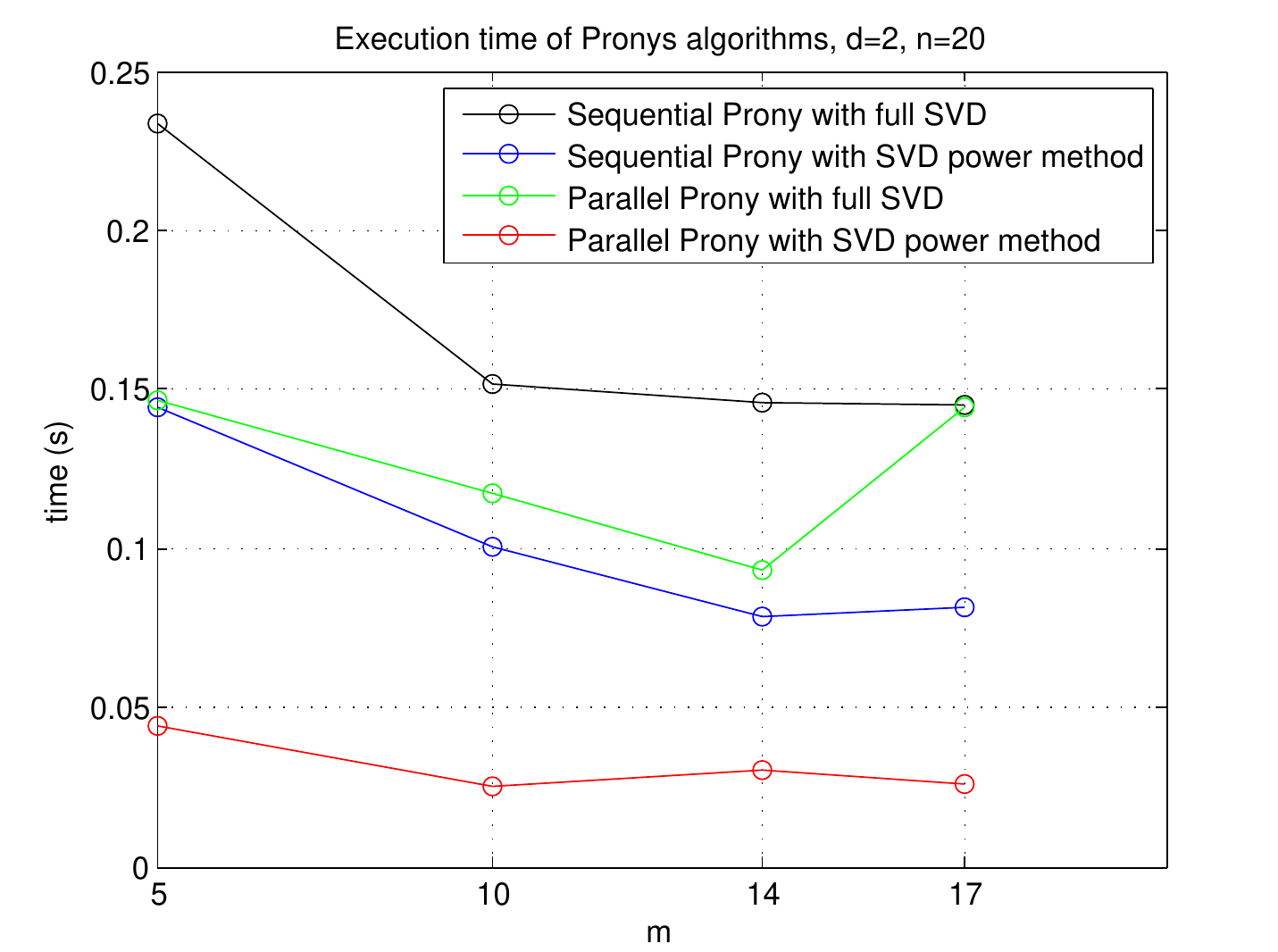} \includegraphics[width=0.47\textwidth, trim=100 245 100 240, clip]{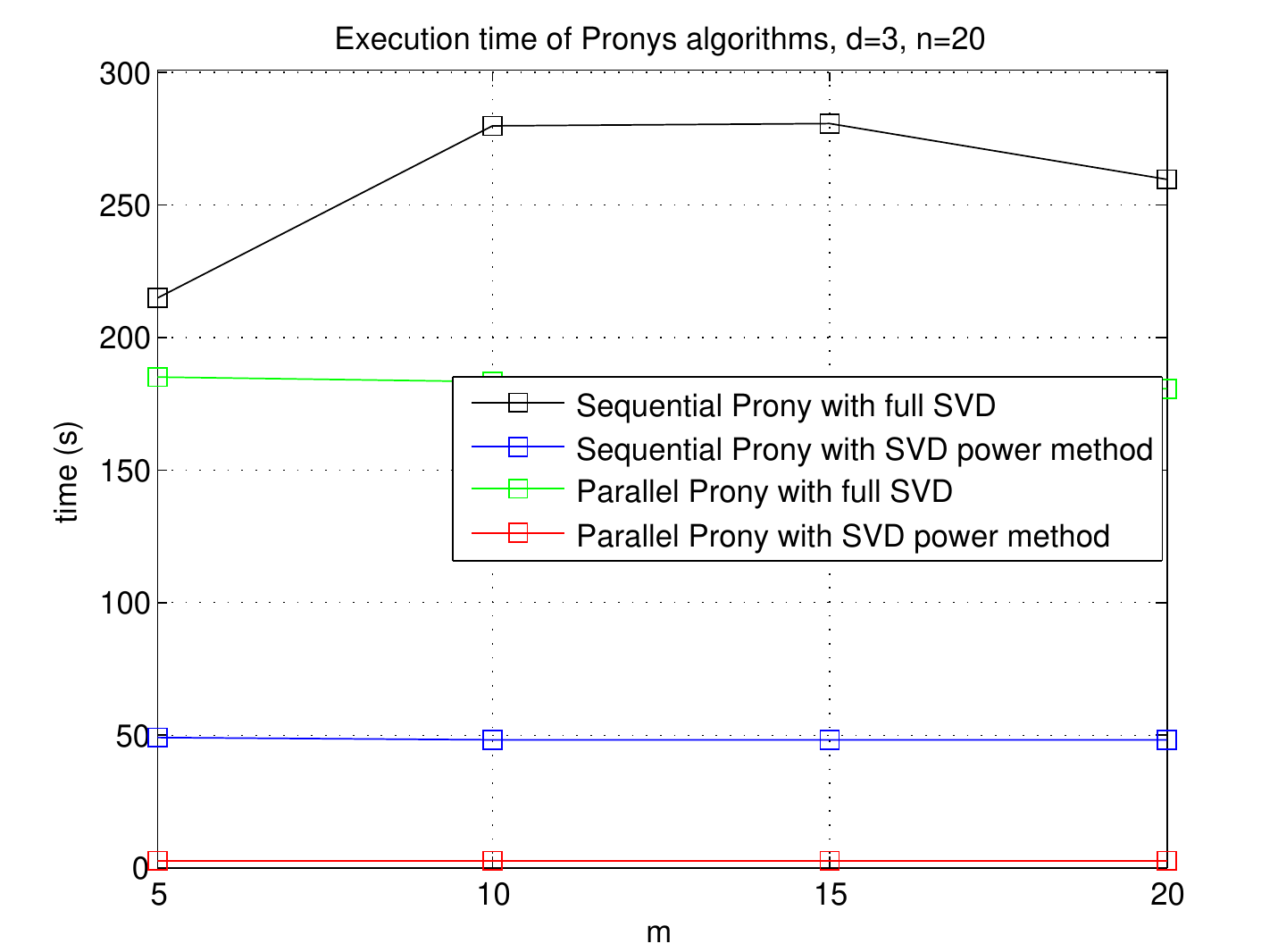}
  \caption{Execution time of four variants of the Prony's method, for $d=2$ and $d=3$.}\label{fig_time_prony}
\end{figure}

\begin{figure}[htbp]
  \centering
  \includegraphics[width=0.8\textwidth, trim=100 245 100 240, clip]{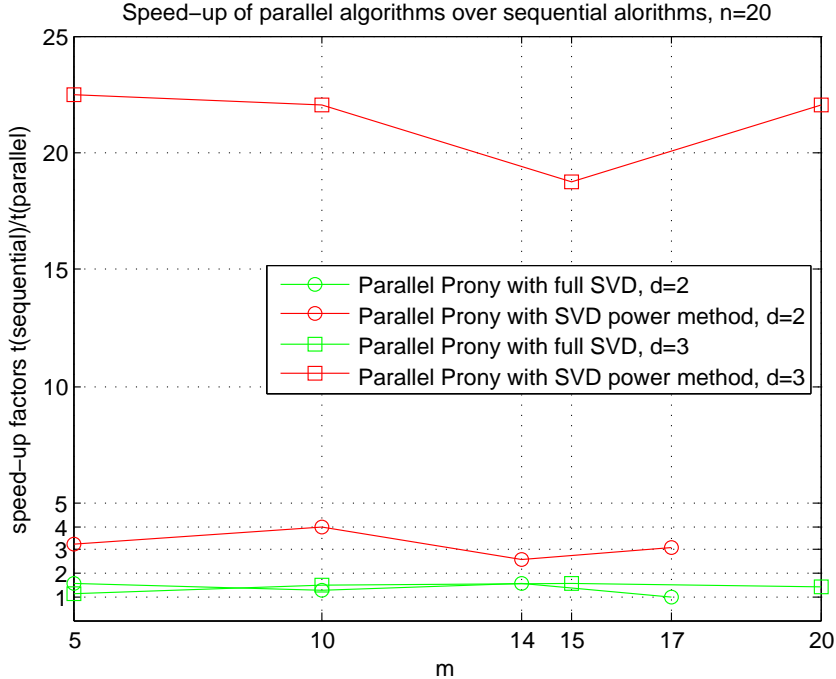}
  \caption{Speed-up factors of parallel over sequential versions of the Prony's method, for $d=2$ and $d=3$.}\label{fig_speedup_prony}
\end{figure}

Next, we are going to compare four variants of the Prony's method
\begin{enumerate}
\item sequential variant with full SVD
\item sequential variant with block power method for SVD
\item parallel variant described in Figure \ref{fig_parallel_prony} with full SVD
\item parallel variant described in Figure \ref{fig_parallel_prony} with block power method for SVD.
\end{enumerate}
Sequential algorithms mentioned above execute their subtasks in sequential order, but not all operations are fully sequential. Matrix operations are implicitly parallelized by calls to multithreaded BLAS routines on our machine. The execution times are presented in Figure \ref{fig_time_prony}, and speed-up factors of parallel over sequential versions can be found in Figure \ref{fig_speedup_prony}.
The first obvious conclusion is that there isn't much benefits from the parallel Prony's mehod with full SVD, its execution time compared to execution time of sequential Prony's algorithm was not reduced much (the largest speed-up factor is about 1.6). The reason for this is that SVD algorithm occupy the largest fraction of the total execution time, for both sequential and parallel algorithm. On the other hand, replacing full SVD with the SVD power method in the sequential algorithm produced speed-up factors up to 5.8 for $d=3$, and parallel variant with SVD power method is in some cases over 120 times faster than the original method. This parallel algorithm is much more efficient even than its sequential version, and for $d=3$ it is about 22 times faster. In the sequential Prony's algorithm with SVD power method, the most demanding task is generation of matrices $T$ and $T_{\ell}$, and the SVD algorithm is very fast. In the parallel version generation of matrices elements is performed in parallel, and thus the larger speed-up factor is obtained.

In the next test round we generate the sampling values from a predetermined exponential sum again but we also introduced noise. The data were taken for $d=3$, $n=20$, and $m=5$, as
$$f(k)(1+\delta_{k}),$$
where the relative perturbation $\delta_{k}$ is bounded by some tolerance $\varepsilon$. The criteria for the rank determination in this case was the first $i$ that satisfies $\sigma_{i}<\text{tol}\cdot \sigma_{1}$, where in the most cases $\text{tol}=\varepsilon$. We measured relative residual norm for the least squares problem in line 7 of Algorithm \ref{alg_prony}, maximal error in the components of vectors $t_{j}$ for $j=1,\ldots ,m$, and the relative error for the coefficient vector $c$. The results are displayed in Table \ref{tab_accuracy}.
\begin{table}[ht]
\begin{center}
\begin{tabular}{|r|r|r|r|r|}
\hline
$\boldsymbol{\varepsilon}$&\textbf{tol}&$\mathbf{\| \tilde{A}^{T}\tilde{c}-\tilde{f}\|_{2}/\| \tilde{f}\|_{2}}$&$\mathbf{\boldsymbol{\max}_{j=1,\ldots ,d \atop i=1,\ldots ,d}|\tilde{t}_{j}(i)-t_{j}(i)|}$&$\mathbf{\| \tilde{c}-c\|_{2}/\| c\|_{2}}$\\
\hline
0&$N\varepsilon_{M}$&1.40484e-14&4.38538e-15&7.67293e-13\\
\hline
1e-9&1e-9&3.00100e-10&1.13784e-11&9.50551e-10\\
\hline
1e-6&1e-6&3.00100e-07&1.13789e-08&9.50556e-07\\
\hline
1e-3&1e-4&2.99893e-04&1.13424e-05&9.52641e-04\\
\hline
\end{tabular}
\end{center}
\caption{Accuracy results for the Prony's method in case when $d=3$, $n=20$, and $m=5$.}\label{tab_accuracy}
\end{table}

As we can see, errors in parameters $t_{j}$ and coefficients $c_{j}$ are proportional to the error in the sample values, as numerical analysis in section \ref{sec_num_anal} suggested. Only for $\varepsilon=10^{-3}$ the tolerance for rank determination has to be smaller, because otherwise the algorithm detected smaller numerical rank $m=4$ and the errors were large. The relative residual norm presented in Table \ref{tab_accuracy} can be used as operational accuracy estimate within the algorithm, implying whether the rank is well determined or not.

\section{Conclusion}
In this paper we proposed a parallel version of Prony's method with multivariate matrix pencil approach, which with appropriate choice of the SVD algorithm for the matrix $T$ (\ref{eq_T_svd}) turned out to be up to 120 times faster than the original algorithm on our computing platform. Since we expect matrix $T$ to have low rank, the most convenient SVD algorithms are the Lanczos bidiagonalization method and the block power method for SVD. The subtasks of the Prony's method are distributed over CPU threads and GPU block of threads, where the CPU is dedicated to more complex tasks such as computing SVD, and GPU to more massive tasks. We also provided a detailed numerical analysis of the method's performance in case of noisy input data, showing that the forward error is proportional to the error in input data, and this was confirmed by numerical tests.

%Spomenuti Lanczosevu bidijagonalizaciju i potrebnu reortogonalizaciju kad ne znamo rang --> istestirati i vjerojatno je sporije.
%Onda opisati SVD metodu potencija kad znamo otprilike rang --> vjerojatno brža.

\end{document}